\documentclass[11pt,reqno]{amsart}
\usepackage[utf8]{inputenc}

\usepackage[all]{xy}
\usepackage{amscd}
\usepackage{mathrsfs}
\usepackage{amssymb}
\usepackage{amsthm}
\usepackage{amsmath}
\usepackage{enumerate}
\usepackage{indentfirst}
\usepackage{dsfont}
\usepackage{tikz-cd}
\usepackage{mathtools}
\usepackage{mathdots}

\usepackage[active]{srcltx}
\usepackage{xcolor}
\usepackage{comment}

\usepackage{hyperref}

\usepackage{xypic}

\theoremstyle{plain}
\newtheorem{thm}{Theorem}[section]
\newtheorem{prop}[thm]{Proposition}
\newtheorem{lem}[thm]{Lemma}
\newtheorem{cor}[thm]{Corollary}
\theoremstyle{definition}
\newtheorem{defn}[thm]{Definition}
\newtheorem{ex}[thm]{Example}
\newtheorem{rem}[thm]{Remark}

\newtheorem{conv}[thm]{Convention}

\newcommand{\fka}{{\mathfrak a}}

\newcommand{\Z}{\mathbb{Z}}
\newcommand{\D}{\mathbb{D}}
\newcommand{\Q}{\mathbb{Q}}

\newcommand{\A}{\mathbb{A}}
\newcommand{\C}{\mathbb{C}}
\newcommand{\R}{\mathbb{R}}

\newcommand{\cA}{\mathcal{A}}
\newcommand{\cE}{\mathcal{E}}
\newcommand{\cO}{\mathcal{O}}

\newcommand{\cP}{\mathcal{P}}

\newcommand{\cI}{\mathcal{I}}
\renewcommand{\cL}{\mathcal{L}}

\newcommand{\Gal}{\operatorname{Gal}}

\newcommand{\Frob}{\tu{Frob}}
\newcommand{\Std}{\tu{Std}}

\newcommand{\sgn}{\tu{sgn}}
\newcommand{\tr}{\tu{tr}\,}
\newcommand{\Hom}{\tu{Hom}}

\newcommand{\End}{\tu{End}}
\newcommand{\Aut}{\tu{Aut}}

\newcommand{\disc}{\tu{disc}}

\newcommand{\temp}{\tu{temp}}
\newcommand{\cusp}{\tu{cusp}}

\newcommand{\SL}{\tu{SL}}

\newcommand{\GL}{\tu{GL}}

\newcommand{\Lie}{\tu{Lie}\,}

\newcommand{\triv}{{\mathbf{1}}}

\newcommand{\Ad}{\tu{Ad}}
\newcommand{\ad}{\tu{ad}}

\newcommand{\Int}{\tu{Int}}

\newcommand{\Irr}{\tu{Irr}}

\newcommand{\vol}{\tu{vol}}

\newcommand{\Res}{\textup{Res}}

\newcommand{\tu}{\textup}
\newcommand{\ol}{\overline}
\newcommand{\ra}{\rightarrow}
\newcommand{\isom}{\stackrel{\sim}{\ra}}

\renewcommand{\hat}{\widehat}

\numberwithin{equation}{subsection}

\begin{document}

\title{The local intertwining relation for inner forms of $\GL_N$}

\author[
H. Atobe, T. Kaletha, A. M{\'i}nguez, S. W. Shin
]
{Hiraku Atobe, 
Tasho Kaletha, 
Alberto M\'inguez
\and
Sug Woo Shin
}

\address{
Department of Mathematics, Kyoto University, 
Kitashirakawa-Oiwake-cho, Sakyo-ku, Kyoto, 606-8502, Japan
}
\email{
atobe@math.kyoto-u.ac.jp
}

\address{
Bonn, Germany
}
\email{
kaletha@math.uni-bonn.de 
}

\address{
Departamento de \'Algebra and Instituto de Matem\'aticas (IMUS), Universidad de Sevilla, C/ Tarifa s/n, 41012 Sevilla, Espa\~{n}a
}
\email{
aminguez@us.es
}

\address{
Department of Mathematics, UC Berkeley, 
Berkeley, CA 94720, USA 
/ 
Korea Institute for Advanced Study, 
Seoul 02455, Republic of Korea
}
\email{
sug.woo.shin@berkeley.edu
}

\date{\today}

\begin{abstract}
The local intertwining relation for inner forms of $\GL_N$ over a local field asserts that the normalized intertwining operator for every parabolically induced representation arising from a local $A$-parameter is the identity operator. Our main theorem verifies the claim over local fields of characteristic zero. This will serve as a vital ingredient for completing the local intertwining relation and endoscopic classification for non-quasi-split unitary groups.
\end{abstract}

\maketitle

\setcounter{tocdepth}{1}
\tableofcontents

\section{Introduction}

The local intertwining relation (LIR) is a key identity in representation theory of reductive groups over local fields. In Arthur’s formulation of endoscopic classification as in \cite{Arthur}, it is one of the two main local theorems, the other being the local classification theorem by packets satisfying the endoscopic character relation (ECR) and other properties. In \emph{op.~cit.}, LIR facilitates the inductive proof of the classification theorem when dealing with parabolically induced representations. However, it is worth stressing that LIR is not a consequence of either the classification theorem or ECR but rather provides further information. 

In the fundamental case of general linear groups, the classification theorem amounts to the local Langlands correspondence \cite{Hen00,HT01,LRS}, whereas LIR asserts that certain normalized intertwining operators are the identity operators. To prove the latter statement, one uses Whittaker models in the tempered case, and deduces the general case via meromorphic continuation \cite{Arthur,AGIKMS}.

The goal of this paper is to prove LIR for inner forms of $\GL_N$ over local fields of characteristic zero. (See the last paragraph of \S\ref{ss:literature} for the characteristic zero assumption.)
The classification theorem in this case amounts to the well-known (local) Jacquet--Langlands correspondence \cite{Rog83,DKV84} (also see \cite{BadRen}), the simplest instance of endoscopic Langlands functoriality. Again, LIR is a statement built on top of the Jacquet--Langlands correspondence, that some intertwining operators for inner forms of $\GL_N$ are the identity operators under a suitable normalization. We are going to verify this by combining both local and global techniques, detailed below. 

Our paper has two main motivations, apart from the intrinsic reason that inner forms of $\GL_N$ arise naturally in many problems. First, we need LIR for inner forms of $\GL_N$ to complete the endoscopic classification for (non-pure) inner forms of quasi-split unitary groups initiated in \cite{KMSW}, because such inner forms of $\GL_N$ appear as local components of adelic unitary groups. Second, the current paper is meant to be a gentle introduction to LIR by showcasing various ideas and techniques, since most ingredients for proving LIR for classical groups appear here but in a simpler form.

\subsection{Main result} 
Let $D$ be a central division algebra over a local field $F$ of characteristic zero with $[D:F]=d^2$. Consider $G=\GL_m(D)$ as an inner form of $G^*=\GL_{dm}$ over $F$. Let $P$ be a proper parabolic subgroup of $G$ with Levi factor $M$. Let $\psi_M$ be a local Arthur parameter of $M$, which is assumed to be $D$-compatible (Definition \ref{def:D-compatible}) so that it determines an irreducible unitary representation $\pi_M$ of $M(F)$ via the local Langlands correspondence for $M$ (extended to Arthur parameters). Write $\cI_P(\pi_M)$ for the normalized parabolic induction; this unitary representation of $G(F)$ is known to be irreducible (since $G$ is an inner form of a general linear group). See \cite{Ber84,Sec09} and \cite{VoganUnitaryDual}. Let $N(M,G)$ denote the normalizer of $M$ in $G$. Write $W_{\psi_M}$ for the subgroup of elements $w\in N(M,G)/M$ such that the natural action of $w$ on $\psi_M$ fixes the isomorphism class of $\psi_M$. We can define the normalized intertwining operator
\[
R_P(w,\psi_M)\in \End_{G(F)}(\cI_P(\pi_M))
\]
by putting together three ingredients (see \S\ref{ss:normalized-intertwining-operators} below for details): (i) meromorphic continuation of the intertwining integral, (ii) a normalizing factor consisting of values of the adjoint $L$ and $\epsilon$-functions arising from $\psi_M$, (iii) an $M$-equivariant linear operator between $\pi_M$ and its $w$-twist constructed by twisted endoscopy. 
Unlike the case of general reductive groups, $R_P(w,\psi_M)$ is a scalar operator since $\cI_P(\pi_M)$ is always irreducible. However, it is an important and subtle problem to identify the scalar. 
Here is our main result.

\begin{thm}
For every $w\in W_{\psi_M}$, the operator $R_P(w,\psi_M)$ is the identity.
\end{thm}

Our theorem is equivalent to the local intertwining relation for
unitary groups in \cite[\S2.6]{KMSW} when the quadratic algebra is
split. The formulation there is parallel to Arthur's
\cite[Thm.~2.4.1]{Arthur}. We refer to \cite[\S\S2.6--2.7]{KMSW}
for the comparison of normalizations. It is also
helpful to think of our theorem as analogous to \cite[Cor.~6.4.5, 7.3.5]{Arthur}
(corollaries of Arthur's LIR).

\subsection{Sketch of proof} 

By multiplicativity and compatibility with parabolic induction, we reduce the general case to proving LIR for $G=\GL_{m}(D)$ with $m$ even, $M=\GL_{m/2}(D)\times\GL_{m/2}(D)$, $w$ exchanging the two factors, and $\psi_M=(\eta,\eta)$ for an
irreducible $D$-compatible $A$-parameter $\eta$ of $\GL_{m/2}(D)$.

Further progress is possible by utilizing two global ingredients proven by the trace formula: the global intertwining relation and a globalization of an irreducible representation over a local field to an automorphic representation. The basic idea is that, in the global setting, if we know LIR at all but one place then we can deduce LIR at the remaining place via the global intertwining relation. With the additional weapons at our disposal, we can choose our favorite battlefield, that is, we can further reduce to the case where $M=\GL_1(D)\times \GL_1(D)$ (where $D$ is still a central division algebra) over a non-Archimedean~$F$, and $\psi_M$ has a very specific form.

We need to face the problem head-on at this point, first dealing with tempered $\psi_M$ and then non-tempered $\psi_M$. 
The tempered case reduces to the ``Steinberg'' parameter $\psi_M$, which corresponds to the trivial representation of $M(F)$. Then the computation of $R_P(w,\psi_M)$ is essentially the same as the analogous computation for principal series of $\GL_2(F)$. The non-tempered case is the most subtle part of this paper, where the correct formulation of LIR requires a tricky twisted-endoscopic sign, presenting a new difficulty. We overcome it by reducing to the case of co-tempered parameters by globalization and employing the Aubert--Zelevinsky (AZ) duality to relate to LIR for tempered parameters; here we exploit the fact that co-tempered parameters become tempered under AZ duality, while irreducible representations of $M(F)=D^\times\times D^\times$ remain invariant. Also unchanged under the duality are the unnormalized intertwining operators, so the proof boils down to comparing the normalizing factors for $\psi_M$ and its AZ dual. We carry out this computation in Lemmas \ref{lem:a-psi} and \ref{lem:rP-ratio}.

The idea to apply AZ duality in the non-tempered LIR originates from \cite[Ch.7]{Arthur}; see also \cite{AGIKMS}. Our strategy illustrates the idea in a simpler setting. On the other hand, the presence of twisted-endoscopic signs, unlike \cite{Arthur,AGIKMS} concerning quasi-split groups, is a new phenomenon due to the non-existence of Whittaker normalization.

\subsection{Organization}

In \S\ref{s:JL} we set up basic notation and recall relevant facts on the local and global Jacquet--Langlands correspondence for inner forms of $\GL_N$. We state the local intertwining relation (LIR) and perform initial reduction steps in \S\ref{s:LIR-statement}. To go further, we bring in global ingredients in \S\ref{s:globalization}. The remaining sections \S\ref{s:proof-tempered} and \S\ref{s:proof-non-tempered} are devoted to completing the proof of LIR in the tempered and non-tempered cases, respectively. In Appendix \ref{akms-reg:appendix} we prove the regularity, invertibility and multiplicativity of the Arthur-normalized intertwining operator for local $A$-parameters.
Appendix~\ref{akms-ind:appendix} proves the compatibility with
parabolic induction used in Lemma~\ref{lem:reduction-to-maximal-Levi}.

\subsection{Relationship to the literature}\label{ss:literature}
We are citing \cite{KMSW}, which has remained a preprint for quite a long time, so we clarify its status and the dependence of our paper on it. We refer to that work largely in two ways: for the local Jacquet--Langlands correspondence (\S\ref{ss:localJL}), and for the definition and basic properties of the normalized intertwining operators (\S\ref{s:LIR-statement}). The former is just a reformulation of the results published by others. The latter extends a construction in \cite[Ch.2]{Arthur} to non-quasi-split inner forms in a purely local setting. The construction depends on the endoscopic classification for proper Levi subgroups, which is given by the local Jacquet--Langlands correspondence in the case at hand; we can entirely avoid considering local unitary groups.

The paper \cite{KMSW} has not been submitted for publication mainly because the authors were waiting for the endoscopic classification of quasi-split unitary groups \cite{Mok} to become unconditional, so that \cite{KMSW} would also become unconditional. Now that a recent preprint \cite{AGIKMS} resolved a key issue, \cite{KMSW} is actively being revised with some updates for publication. Another reason for the delayed publication was to make sure that the non-tempered case for non-quasi-split unitary groups is within sight. This problem amounts to ``[A25]'' (see \cite{AGIKMS}) for such groups. We are in the process of adapting the strategy of \cite{AGIKMS} for quasi-split classical groups to non-quasi-split classical groups. In summary, the submission of \cite{KMSW} and its sequel is on the horizon.

The use of the trace formula in this paper is limited to general linear groups and their inner forms, where \cite{AC89} provides all we need. Unlike in \cite{Arthur,Mok,KMSW}, neither the weighted fundamental lemma nor the stabilization of the general trace formula is required. 

We are limiting ourselves to local fields of characteristic zero because not all of our ingredients are available when the characteristic is positive. (E.g., compare \cite{BadRoc} with \cite{Bad,BadRen}.) That said, we think that LIR is within reach in that case, just requiring further work.

\subsection{Acknowledgments}
HA was partially supported by JSPS KAKENHI Grant Number 23K12946. AM was partially funded by the Principal Investigator project PAT-4832423 of the Austrian Science Fund (FWF) and the project ATR2024-154613 funded by MICIU/AEI/10.13039/501100011033.
SWS was partially supported by NSF grant DMS-2401353, NSF RTG grant DMS-2342225, and Simons Travel Grant. We thank Johannes Droschl for useful communications. 
After completing the draft without AI (before writing the appendix), we used ChatGPT for proofreading and soliciting suggestions for minor corrections and improvements. ChatGPT also assisted us with improving the proof of Lemma \ref{lem:reduction-to-maximal-Levi}, completing the proofs in Appendix~\ref{akms-reg:appendix} on basic properties of normalized intertwining operators, as well as writing Appendix~\ref{akms-ind:appendix} to give a self-contained justification of a fact in \cite[\S2.7]{KMSW}.

\subsection{Notation}\label{ss:notation}

When $F$ is a local field, write $\Gamma_F$ (resp.~$W_F$) for the full Galois (resp.~Weil) group over $F$. Often we write $\Gamma$ for $\Gamma_F$.
Denote by $|\cdot| \colon F^\times\ra \R_{>0}$ the usual absolute value. 
If $F$ is $p$-adic, it sends a uniformizer to $q^{-1}$, 
where $q$ is the cardinality of the residue field of $F$. We keep writing $|\cdot|$ for the absolute value character on $W_F$ via local class field theory (sending a uniformizer to a geometric Frobenius element). 
Fix a nontrivial unitary additive character $\psi_F
\colon F\to \C^\times$. 
We often suppress $\psi_F$ from the notation.

We usually write $D$ for a division algebra with center $F$, and $d:=[D:F]^{1/2}$. For the reductive group $G=\GL_m(D)$, we write $G^*$ for the quasi-split inner form $\GL_{dm}$. 
Composing with determinant and/or reduced norm, $|\cdot|$ is extended to a character on $\GL_m(F)$ and more generally on $\GL_m(D)$ for any $m\in\Z_{\ge1}$. 
Once we fix an isomorphism $\xi \colon M_d(\overline{F}) \xrightarrow{\sim} D \otimes_{F} \overline{F}$, 
we extend it to an algebraic isomorphism 
$\xi \colon M_{dm}(\overline{F}) \xrightarrow{\sim} M_m(D \otimes_{F} \overline{F})$
which gives rise to an $\overline{F}$-isomorphism
\[
\xi \colon G^*_{\overline{F}}  \xrightarrow{\sim} G_{\overline{F}}.
\]
We fix a natural splitting $\tu{spl}_{G^*} = (B^*, T^*, \{X_i\}_{1 \leq i \leq dm-1})$, 
where
$B^*$ is the upper-triangular Borel subgroup of $G^*$, 
$T^*$ is the diagonal torus 
and $X_i \in \tu{Lie}(G^*)$ is the matrix whose $(i,i+1)$-component is $1$ and the other components are $0$.

Let $\triv$ denote the trivial representation of a group (which will be clear from the context). By $S_n$ we mean an irreducible algebraic $n$-dimensional representation of $\SL_2(\C)$.

\section{The Jacquet--Langlands correspondence for inner twists of $\GL_N$}\label{s:JL}

\subsection{Parameters}\label{ss:parameters}

Let $G$ be a connected reductive group over a local field $F$. Set $\mathcal L_F:=W_F$ if $F$ is Archimedean and $\mathcal L_F:=W_F\times \SL_2(\C)$ otherwise. We still write $|\cdot|\colon \mathcal L_F\to \C^\times$ for the norm character through the projection onto $W_F$.
We write 
\[
\Irr(G)\supset\Irr_u(G)\supset\Irr_{\tu{temp}}(G)\supset \Irr_2(G)
\]
for the sets of isomorphism classes of irreducible admissible representations of $G(F)$, respectively all, unitarizable, tempered, and square-integrable modulo the center. 

Denote by $\Phi(G)$ (resp.~$\Psi(G)$) the set of isomorphism classes (given by $\hat G$-conjugacy) of $L$-parameters $\mathcal L_F\to {}^L G$ (resp.~$A$-parameters $\mathcal L_F\times \SL_2(\C)\to {}^L G$). For each $\phi\in \Phi(G)$ and $\psi\in \Psi(G)$, define $S_\phi,S_\psi\subset \hat G$ to be the centralizer groups of $\phi,\psi$, respectively.
We have the subset of tempered parameters $\Phi_{\tu{temp}}(G)$ consisting of $\phi\in \Phi(G)$ whose restriction to $W_F$ is represented by a 1-cocycle of $W_F$ valued in a bounded subgroup of $\hat G$. Then $\Phi_{\tu{temp}}(G)$ is identified with the subset of $\Psi(G)$ such that $\psi\in\Psi(G)$ has trivial restriction to the Arthur $\SL_2(\C)$. Define $\Phi_2(G)\subset \Phi_{\tu{temp}}(G)$ and $\Psi_2(G)\subset \Psi(G)$ by the condition that $S_\phi/Z(\hat G)^{\Gamma}$ and $S_\psi/Z(\hat G)^{\Gamma}$ are finite. (In our main setting, the $\Gamma$-action on $\hat G$ will always be trivial, and $S_\phi/Z(\hat G)^{\Gamma}$ and $S_\psi/Z(\hat G)^{\Gamma}$ will be trivial as soon as they are finite.) 
On the other hand, we obtain an enlargement $\Psi^+(G)\supset \Psi(G)$ by dropping the requirement that the restriction of $\mathcal L_F\times \SL_2(\C)\to {}^L G$ to $W_F$ should be a bounded cocycle. We obtain a map
\[
\Psi^+(G)\to \Phi(G),\qquad \psi\mapsto \phi_\psi
\]
by pulling back along the map $\mathcal L_F\to \mathcal L_F\times \SL_2(\C)$, $w\mapsto \left(w,\begin{psmallmatrix}
    |w|^{1/2} & 0\\ 0 & |w|^{-1/2}
\end{psmallmatrix}\right)$.

Now we consider the group $G^*=\GL_N$. The local Langlands correspondence \cite{HT01,Hen00,Lan89} maps $\Irr(G^*)$ bijectively onto $\Phi(G^*)$, inducing bijections $\Irr_{\temp}(G^*)\isom \Phi_{\temp}(G^*)$ and $\Irr_2(G^*)\isom \Phi_2(G^*)$. Let us write $\pi^*_\phi\in \Irr(G^*)$ for the preimage of $\phi\in \Phi(G^*)$. The subset $\Phi_2(G^*)$ of $\Phi(G^*)$ is cut out by the condition that $\phi\in \Phi(G^*)$ is irreducible as an $N$-dimensional representation of $\cL_F$, with bounded image when restricted to $W_F$. 
Each $\psi\in \Psi_2(G^*)$ has the form $\psi=\phi\boxtimes S_b$ for $\phi\in \Phi_2(\GL_m)$ such that $mb=N$; we assign to it the Speh representation $\pi^*_\psi\in \Irr_u(G^*)$ which is the unique (irreducible) Langlands quotient of the normalized induction of $\pi^*_{\phi}|\cdot|^{\frac{b-1}{2}}\otimes \pi^*_{\phi}|\cdot|^{\frac{b-3}{2}}\otimes \cdots \otimes \pi^*_{\phi}|\cdot|^{\frac{1-b}{2}}$. Since $\psi\in \Psi(G^*)$ can be written as 
\[
\psi=\oplus_{i\in I} \psi_i,\quad \psi_i\in \Psi_2(\GL_{n_i}),\quad \textstyle\sum_{i\in I} n_i=N
\]
as an $N$-dimensional representation of $\mathcal L_F\times \SL_2(\C)$,
we can define $\pi^*_\psi$ to be the normalized parabolic induction of $\otimes_{i\in I} \pi^*_{\psi_i}$. This defines an injection
\begin{equation}\label{eq:A-packet-GLn}
    \Psi(G^*)\hookrightarrow \Irr_u(G^*),\qquad \psi\mapsto \pi^*_{\psi}.
\end{equation}
To obtain a bijection, consider the intermediate set  $\Psi(G^*)\subset \Psi^+_u(G^*)\subset \Psi^+(G^*)$ consisting of ``generalized'' $A$-parameters $\psi\colon \mathcal L_F\times \SL_2(\C)\to {}^L G^*$ which have the following form:
\begin{equation}\label{eq:psi-decomposition}
\psi=\bigoplus_{i\in I} \psi_i \oplus \bigoplus_{j\in J} \psi_j\otimes \big(|\cdot|^{\epsilon_j}\oplus |\cdot|^{-\epsilon_j}\big),\quad \psi_i\in \Psi_2(\GL_{n_i}), \quad \psi_j\in \Psi_2(\GL_{n_j}), \quad 0<\epsilon_j<\tfrac{1}{2},
\end{equation}
where $\sum_{i\in I}n_i+2\sum_{j\in J}n_j=N$. We define $\pi^*_{\psi}$ to be the normalized parabolic induction of 
\[\otimes_{i\in I}\pi^*_{\psi_i}\otimes (\otimes_{j\in J}\pi^*_{\psi_j}|\cdot|^{\epsilon_j}\otimes \pi^*_{\psi_j}|\cdot|^{-\epsilon_j}).
\]
Then $\psi\mapsto \pi^*_{\psi}$ is a bijection $\Psi^+_u(G^*)\isom \Irr_u(G^*)$ thanks to the theorems by Tadi\'c and Vogan \cite{TadicUnitaryDual,VoganUnitaryDual}. 

When $F$ is non-Archimedean, each $\psi\in \Psi_2(G^*)$ can be uniquely written as
\begin{equation}\label{eq:sq-int-psi}
  \psi=\rho\boxtimes S_a\boxtimes S_b,
\end{equation}
where 
$\rho$ is an irreducible unitary $m$-dimensional representation of $W_F$, 
and $a,b \in \Z_{\ge 1}$ such that $N=mab$.
Here $S_a$ and $S_b$ are representations of $\SL_2(\C)$ inside and outside $\cL_F$, respectively.
We define the (Aubert) dual parameter
\[
\psi\mapsto \hat \psi:=\psi\circ i, \qquad \psi\in \Psi^+_u(G^*),
\]
where $i$ is the automorphism on $\cL_F\times \SL_2(\C)=W_F\times \SL_2(\C)\times \SL_2(\C)$ given by $(w,g_1,g_2)\mapsto (w,g_2,g_1)$. For example, $\hat \psi=\rho\boxtimes S_b\boxtimes S_a$ for $\psi$ of the form \eqref{eq:sq-int-psi}. This involution is compatible with the Aubert--Zelevinsky involution $\pi^*\mapsto \hat \pi^*$ on $\Irr(G^*)$: 
in terms of \eqref{eq:A-packet-GLn}, 
\begin{equation}\label{eq:Apacket-vs-Aubert}
    \hat \pi^*_{\psi} =  \pi^*_{\hat\psi}.
\end{equation}
See \cite{MW86} or \cite{LM14}.
Finally, we remark that all the definitions above extend to a finite product of general linear groups (e.g., Levi subgroups in $\GL_N$) in the obvious manner.

\subsection{Extended pure inner twists}\label{ss:epit}

Let $G^*$ be a connected quasi-split reductive group over a local or global field $F$ of characteristic $0$. 
For a finite Galois extension $K/F$, consider the extension 
\[
1\to \D_{K/F}(K)\to \cE(K/F)\to \Gal(K/F)\to 1
\]
as in \cite[\S2]{Kot14}, where $\D_{K/F}$ is a certain group of multiplicative type over $F$. 
An \emph{extended pure inner twist} $(G,\xi,z)$ of $G^*$ consists of a connected reductive group $G$ over $F$, an $\ol F$-isomorphism $\xi\colon G^*_{\ol F}\to G_{\ol F}$, and a 1-cocycle $z\colon \cE(K/F)\to G^*(K)$ which is algebraic and basic in the sense of \cite{Kot14} for some finite Galois extension $K/F$, cf.~\cite[\S0.3.1]{KMSW}, such that the map $z_{\ad}\colon \Gal(K/F)\to G^*_{\ad}(K)$ arising from $z$ equals $\sigma\mapsto \xi^{-1}\sigma(\xi)$. Two extended pure inner twists $(G_1,\xi_1,z_1)$ and $(G_2,\xi_2,z_2)$ are isomorphic (over $F$) if there exist an $F$-isomorphism $f\colon G_1\isom G_2$ and $g\in G^*(K)$ for a sufficiently large finite Galois extension $K/F$ such that $\xi_2^{-1}\circ f\circ \xi_1=\Int(g)$ and $z_2(\sigma)=g z_1(\sigma) \sigma(g)^{-1}$ for $\sigma\in \cE(K/F)$, where $\sigma$ acts on $g$ through $\Gal(K/F)$. The isomorphism classes of extended pure inner twists of $G^*$ as $K$ varies are classified by $B(F,G^*)_{\tu{basic}}$, defined as the colimit of all algebraic and basic 1-cocycles $\cE(K/F)\to G^*(K)$ modulo coboundaries. Forgetting the 1-cocycle $z$ from the above, we recover the well-known notion of \emph{inner twist} $(G,\xi)$ of $G^*$ over $F$ together with the notion of isomorphisms.

We have natural maps
\begin{equation}\label{eq:B-to-Bad}
 B(F,G^*)_{\tu{basic}}\to B(F,G^*_{\ad})_{\tu{basic}} = H^1(F,G^*_{\ad}),   
\end{equation}
where the latter classifies inner twists of $G^*$ over $F$ up to isomorphism. The arrow is known to be surjective when the center $Z_{G^*}$ is connected. That is, every inner twist can be promoted to an extended pure inner twist (up to isomorphism) in that case.

This paper is mainly concerned with extended pure inner twists for $G^*=\GL_N$. In this case, we can apply the Hasse invariant map to $H^1(F,G^*_{\ad})$ and identify \eqref{eq:B-to-Bad} with
\[
\Z\twoheadrightarrow \tfrac{1}{N}\Z/\Z
\]
if $F$ is non-Archimedean, whereas if $F$ is Archimedean, \eqref{eq:B-to-Bad} becomes
\[
\Z\twoheadrightarrow \begin{cases}
    \tfrac{1}{2}\Z/\Z, & F=\R,~2\mid N,\\
    0, & F=\R,~2\nmid N,~\mbox{or}~F=\C.
\end{cases}
\]
When $F$ is global, there is a natural exact sequence
\[
0 \to B(F,G^*)_{\tu{basic}}\to \bigoplus_{v} B(F_v,G^*)_{\tu{basic}}\stackrel{\sum_v \kappa_v}{\longrightarrow} X^*(Z_{\hat G^*})=\Z\to 0,
\]
where the second arrow is the localization map over all places $v$ of $F$, and $\kappa_v:B(F_v,G^*)_{\tu{basic}}\rightarrow X^*(Z_{\hat G^*})$ is the Kottwitz morphism at each $v$. If $v$ is finite, $\kappa_v$ is an isomorphism; if $v$ is infinite, it is injective.

Now let $(G,\xi,z)$ be an extended pure inner twist of $G^*=\GL_N$ over a local or global field $F$. Consider two parabolic pairs $(M,P)$ of $G$ and $(M^*,P^*)$ of $G^*$ over $F$, which become $G(\ol F)$-conjugate via $\xi$. By \cite[\S0.4.8]{KMSW} (which allows general linear groups, not just unitary groups) we can always modify $(G,\xi,z)$ by an isomorphism such that $\xi(M^*)=M$ and $\xi(P^*)=P$ over $\ol F$. Then $(M,\xi|_{M^*},z)$ is an extended pure inner twist of $M^*$. (In particular $z$ is necessarily valued in $M^*$.)

\subsection{Local Jacquet--Langlands}\label{ss:localJL}

Let $(G,\xi)$ be an inner twist of $G^*=\GL_N$ over a local field $F$ of characteristic $0$. Then $G=\GL_m(D)$ for a central division algebra $D$ over $F$ such that $d:=[D:F]^{1/2}$ satisfies $N=dm$. We have the Kottwitz sign $e(G)=(-1)^{N-m}$.

Let us define a condition on parameters to ensure that the associated $A$-packet for $G$ is nonempty.

\begin{defn}\label{def:D-compatible}
When $F$ is Archimedean (resp.~$p$-adic),
a parameter $\psi\in \Psi^+_u(G^*)$ is said to be \emph{$D$-compatible} if every irreducible constituent $\phi\boxtimes S_b$ (resp.~$\rho\boxtimes S_a\boxtimes S_b$) of the $\cL_F\times \SL_2(\C)$-representation $\psi$ has the property that $b\dim \phi$ (resp.~either $a\dim \rho$ or $b\dim \rho$) is divisible by $d$. 
Define $\Psi^+_u(G^*)_{D\tu{-comp}}$ to be the subset of $D$-compatible parameters in $\Psi^+_u(G^*)$.
\end{defn}

\begin{conv}\label{convention:D-compatible}
When $G=\GL_{m'}(D')$ for $D'=M_{m/m'}(D)$ with $m'\in \Z_{\ge1}$ dividing $m$, we say that $\psi$ is $D'$-compatible to mean it is $D$-compatible, as this is going to be sometimes convenient.
\end{conv}

\begin{ex}
Suppose that $F$ is a $p$-adic field, $d=8$, $G=\GL_1(D)$, $G^*=\GL_8$, and that
$\rho$ is an irreducible unitary $2$-dimensional representation of $W_F$.
The parameters $\textbf{1}\boxtimes S_1 \boxtimes S_8$ and $\rho\boxtimes S_4\boxtimes S_1$ are $D$-compatible but $\rho\boxtimes S_2\boxtimes S_2$ is not, even though they all belong to $\Psi_2(G^*)$.
\end{ex}

The following may be viewed as the construction of (singleton) $A$-packets for $(G,\xi)$.

\begin{prop}\label{prop:local-JL}
Let $F$, $G^*$ and $(G,\xi)$ be as above. There is a map
\begin{equation}\label{eq:JL-for-A-parameters}
    \Psi^+_u(G^*)_{D\tu{-comp}}\to \Irr_u(G),\qquad \psi\mapsto \pi_\psi,
\end{equation}
uniquely characterized by the following property: there exists a sign $a_\psi\in \{\pm1\}$ such that the Harish-Chandra characters of $\pi_\psi$ and $\pi^*_\psi$ are related by
\[
\Theta_{\pi_\psi}(g)=e(G)a_\psi\Theta_{\pi^*_\psi}(g^*),\qquad g\in G(F)~\leftrightarrow~g^*\in G^*(F),
\]
where $g\leftrightarrow g^*$ means that they are regular semisimple elements which become conjugate over $\ol F$ via $\xi$.
\end{prop}

\begin{proof}
This is reformulated from the main local results of \cite{Bad,BadRen}.
\end{proof}

\begin{conv}\label{convention:transfer-factor}
For Jacquet--Langlands transfer, we normalize the
Langlands--Shelstad transfer factor to equal one whenever nonzero
and choose compatible Haar measures on inner forms. That is, if $f^*\in \mathcal{H}(G^*)$ is a transfer of $f\in \mathcal{H}(G)$, and if $\gamma\in G(F)$ and $\gamma^*\in G^*(F)$ are regular semisimple with matching conjugacy classes, then we have the orbital integral identity $O_{\gamma}(f)=O_{\gamma^*}(f^*)$. (Of course $O_{\gamma^*}(f^*)=0$ if there exists no $\gamma$ whose conjugacy class matches that of $\gamma^*$.) Thus the character identity in Proposition \ref{prop:local-JL} can be rewritten as 
$\tr \pi_\psi(f)=e(G)a_{\psi} \tr \pi^*_\psi(f^*)$. See~\cite[Thm.~1.6.4]{KMSW}.
\end{conv}

The sign $a_{\psi}$ can be computed explicitly as in \cite[\S1.6.3]{KMSW}. For later use, let us recall the formula when $F$ is $p$-adic. Take a decomposition $\psi=\oplus_i k_i \psi_i$ into mutually non-isomorphic irreducible representations $\psi_i$ with multiplicities $k_i\in \Z_{\ge1}$. Write $\psi_i=\rho_i\boxtimes S_{a_i}\boxtimes S_{b_i}$. We define $a_{\psi}:=\prod_i a_{\psi_i}^{k_i}$, where $a_{\psi_i}\in \{\pm1\}$ is given as follows. If $d\nmid(b_i \dim \rho_i)$ then $a_{\psi_i}:=1$. If $d \mid (b_i \dim \rho_i)$, take the smallest $s\in \Z_{\ge1}$ such that $d\mid (s \dim \rho_i)$; then $a_{\psi_i}:=1$ if $s$ is odd, and $a_{\psi_i}:=(-1)^{a_i b_i/s}$ if $s$ is even.

\begin{lem}\label{lem:a-psi}
In Proposition \ref{prop:local-JL}, consider the case $G=D^\times$, $F$ non-Archimedean, and $\psi=\rho\boxtimes S_a\boxtimes S_b\in \Psi_2(G^*)_{D\tu{-comp}}$ (so $\min\{a,b\} = 1$). Then $a_\psi=1$ if $b=1$, and $a_{\psi}=(-1)^{b-1}$ if $a=1$.
\end{lem}

\begin{proof}
This is straightforward from the preceding formula for $a_{\psi}$.
\end{proof}

\begin{rem}\label{rem:not-inj-surj}
The map \eqref{eq:JL-for-A-parameters} need not be injective nor surjective. Let us consider examples with $F$ non-Archimedean and $d>1$. If $G=D^\times$ then both $\mathbf{1}\boxtimes S_d \boxtimes S_1$ and $\mathbf{1}\boxtimes S_1\boxtimes S_d$ correspond to the trivial representation of $D^\times$ (which is also the Steinberg representation of $D^\times$). 
The map may not be onto $\Irr_u(G)$.
For example, if $G = \GL_2(D)$, 
the normalized parabolic induction of $|\cdot|^s \otimes |\cdot|^{-s}$ is irreducible and unitary for $0 < s < d/2$, 
and it is not in the image of the map \eqref{eq:JL-for-A-parameters} for $1/2 < s < d/2$. Another interesting example is the following. Let $\pi$ denote the Steinberg representation of $\GL_2(D)$. Then the normalized parabolic induction of $\pi|\cdot|^s \otimes \pi|\cdot|^{-s}$ is again irreducible and unitary for $0 < s < d/2$. At $s=1/2$, it is of Arthur type, with associated parameter $\mathbf{1}\boxtimes S_{2d}\boxtimes S_2$. By contrast, the Langlands quotient of the reducible parabolic induction from $\pi|\cdot|^{d/2}\otimes \pi|\cdot|^{-d/2},$ namely the unitary Speh representation ${\rm Speh}(\pi,2)$ of $\GL_4(D)$, is not of Arthur type, since it does not arise as the local Jacquet--Langlands transfer of any unitary representation of $\GL_{4d}(F)$. These non-Archimedean examples follow from the definitions and transfer formulas in \cite[\S\S2.4, 2.9, 3.2--3.4, especially (2.1), (2.2), and (3.9)]{Bad}. In that notation $s(\mathbf{1}_{D^\times})=s(\pi)=d$, $\mathrm{Speh}(\pi,2)=u'(\pi,2)$, whereas $|\operatorname{LJ}|(u(\mathrm{St}_{2d},2))=\bar u(\pi,2)
=\pi|\cdot|^{1/2}\times\pi|\cdot|^{-1/2}$. The factorization formulas also show that $u'(\pi,2)$ is not in the image of $|\operatorname{LJ}|$ on the unitary dual. For an analogous Archimedean complementary-series example, see the bottom paragraph of \cite[p.~1118]{BadRen}.
We have the interpretation 
\begin{equation}\label{eq:apsi-via-spsi}
    a_\psi=\langle s_\psi,\pi_\psi\rangle
\end{equation}
in terms of $s_\psi=\psi(1_{\cL_F},-1)\in S_\psi:=\tu{Cent}(\psi,\hat G^*)$ and the character $\langle \cdot,\pi_\psi\rangle \colon S_\psi\to \C^\times$ as in \cite[Thm.~1.6.1]{KMSW}; see \S1.6.3 therein for explanation.
\end{rem}

\begin{rem}
For $\psi \in \Psi_u^+(G^*)_{D\tu{-comp}}$, 
we define the $A$-packet associated to $\psi$ 
by the singleton $\Pi_\psi = \{\pi_\psi\}$.
When $\psi\in \Psi(G^*)$ is not $D$-compatible, the $A$-packet $\Pi_\psi\subset \Irr_u(G)$ is considered to be empty, in which case the equality in Proposition \ref{prop:local-JL} still holds in that the right hand side vanishes whenever $g^*$ corresponds to some regular semisimple element of $G(F)$. See \cite[Thm.~1.6.4]{KMSW}, part 2. We will not use this fact.
\end{rem}

\subsection{Global Jacquet--Langlands}\label{ss:globalJL}

Let $\dot F$ be a number field. For a reductive group $H$ over $\dot F$,
write $\cA_{\cusp}(H)$ and $\cA_2(H)$ for the sets of its cuspidal and discrete
unitary automorphic representations, respectively. We require their central characters
to be trivial on $A_H(\R)^0$, where $A_H$ is the maximal $\Q$-split torus in
the center of $\Res_{\dot F/\Q}H$.

For $\dot G^*=\GL_{N,\dot F}$, define $\Psi_2(\dot G^*)$ as the set of
formal products
\[
 \dot\psi=\dot\phi\boxtimes S_b,\qquad
 \dot\phi\in\cA_{\cusp}(\GL_{a,\dot F}),\qquad ab=N.
\]
The theorem of M\oe glin--Waldspurger \cite{MW89} gives a bijection
$\dot\psi\mapsto\pi^*_{\dot\psi}$ from $\Psi_2(\dot G^*)$ to
$\cA_2(\dot G^*)$, cf.~\cite[Thm.~1.3.3]{Arthur}. Write $\dot G_v^*$ for
$\GL_{N,\dot F_v}$. For each place $v$ there is a localization map
\begin{equation}\label{eq:localize-A-parameters}
 \Psi_2(\dot G^*)\longrightarrow\Psi_u^+(\dot G_v^*),\qquad
 \dot\phi\boxtimes S_b\longmapsto
 \dot\psi_v=\dot\phi_v\boxtimes S_b,
\end{equation}
where $\dot\phi_v$ is viewed as an $L$-parameter via the local Langlands
correspondence. The representation attached to $\dot\psi_v$ by
\S\ref{ss:parameters} is the local component $\pi^*_{\dot\psi,v}$.
The generalized Ramanujan conjecture would imply
$\dot\psi_v\in\Psi(\dot G_v^*)$ at every place.

Let $(\dot G,\dot\xi)$ be an inner twist of $\dot G^*$ over $\dot F$. We have $\dot G=\GL_m(\dot D)$ for a central division algebra $\dot D$ over $\dot F$ with $d:=[\dot D:\dot F]^{1/2}=N/m$. Write $\dot G_v:=\dot G_{\dot F_v}$ and $\dot D_v:=\dot D\otimes_{\dot F}\dot F_v$. We say $\dot\psi\in \Psi_2(\dot G^*)$ is \emph{$\dot D$-compatible} if $\dot\psi_v$ is $\dot D_v$-compatible at every place $v$ of $\dot F$; we follow Convention \ref{convention:D-compatible} as $\dot D_v$ need not be a division algebra. (The condition is vacuous at $v$ if $\dot D$ splits at $v$.)
Denote by $\Psi_2(\dot G^*)_{\dot D\tu{-comp}}$ the subset of $\dot D$-compatible parameters in $\Psi_2(\dot G^*)$. Define
\[
\pi_{\dot \psi}:=\otimes'_v \pi_{\dot \psi_v}\qquad \mbox{for}~\dot\psi\in \Psi_2(\dot G^*)_{\dot D\tu{-comp}}.
\]

\begin{prop}\label{prop:global-JL}
The map $\dot\psi\mapsto\pi_{\dot\psi}$ is a bijection from
$\Psi_2(\dot G^*)_{\dot D\tu{-comp}}$ to $\cA_2(\dot G)$.
Each $\pi_{\dot\psi}$ occurs with multiplicity one and transfers to
$\pi^*_{\dot\psi}$ under global Jacquet--Langlands.
\end{prop}

\begin{proof}
This is a reformulation of the main global theorems of \cite{Bad,BadRen}.
\end{proof}

\section{The local intertwining relation}\label{s:LIR-statement}

\subsection{Normalized intertwining operators}\label{ss:normalized-intertwining-operators}

We are going to recall the definition of normalized intertwining operators following \cite{Arthur,KMSW} in the case of inner forms of general linear groups. Certain things simplify in this case. First, the Keys--Shahidi $\lambda$-factors are always equal to $1$ (thus they will be omitted from the definition). Second, the parabolic induction of an irreducible unitary representation is always irreducible; hence a self-intertwining operator is always a scalar.

We start with the following data:
\begin{itemize}
    \item $F$ is a local field of characteristic $0$,
    \item $\psi_F\colon F\to \C^\times$ is a nontrivial unitary additive character,
    \item $(G,\xi,z)$ is an extended pure inner twist of $G^*:=\GL_N$ over $F$; then $G=\GL_m(D)$ for a central division algebra $D$ over $F$ such that $d:=[D:F]^{1/2}=N/m$,
    \item $(M,P)$ is a proper parabolic pair of $G$ over $F$,
    \item $(M^*,P^*)$ is the unique standard parabolic pair of $G^*$ conjugate to $\xi^{-1}(M,P)$ over $\ol F$; it is proper and defined over $F$ (if we write $M=\prod_i \GL_{m_i}(D)$ then $M^*=\prod_i \GL_{dm_i}$ as a block diagonal subgroup of $G^*$),
    \item $(\hat M^*,\hat P^*)$ is the corresponding standard parabolic pair of $\hat G^*$,
    \item $N,N^*,\hat N^*$ are the unipotent radicals of $P,P^*,\hat P^*$,
    \item $A_M$ is the maximal $F$-split torus in the center of $M$,
    \item $A_{\hat M^*}$ is the maximal torus in $Z_{\hat M^*}$ with trivial $\Gamma$-action (so $A_{\hat M^*}=Z_{\hat M^*}$ in our case),
    \item $X^*_F(M)$ is the group of $F$-rational characters of $M$,
    \item $\fka_M^*:=X^*_F(M)\otimes_{\Z}\R=X^*_F(M^*)\otimes_{\Z}\R$ and $\fka_{M,\C}^*:=\fka_M^*\otimes_{\R}\C$,
    \item $\chi_\lambda\colon M(F)\to \C^\times$ and $\zeta_\lambda\colon W_F\to Z(\hat M^*)$ are continuous morphisms determined by $\lambda\in \fka_{M,\C}^*$,
    \item $\psi_M\in \Psi(M^*)$ and $\psi_{M,\lambda}:=\psi_M\otimes \zeta_\lambda\in \Psi^+(M^*)$,
    \item $\psi\in \Psi(G^*)$ is the image of $\psi_M$, 
        \item $S_\psi^{\tu{rad}}:=(S_\psi\cap \hat{G}^*_{\tu{der}})^\circ$, where $\hat{G}^*_{\tu{der}}$ is the derived subgroup of $\hat G^*$. 
    \end{itemize}
    
    From now on, we assume that $\psi_M$ is $D$-compatible.
    
    \begin{itemize}
    \item $\pi^*_M = \pi^*_{\psi_M}\in  \Irr_u(M^*(F))$ denotes the representation of $M^*(F)$ associated to $\psi_M$,
    \item 
    write $\pi_M = \pi_{\psi_M} \in \Irr_u(M(F))$,
    \item $\pi_{M,\lambda}:=\pi_M\otimes \chi_\lambda$, where $\lambda\in \fka_{M,\C}^*$,
    \item $\mathcal{I}_P(\pi_{M,\lambda})$ for the normalized parabolic induction of $\pi_{M,\lambda}$ from $P$ to $G$,
    \item $Z(H,K),~N(H,K)$ denote the centralizer and normalizer of a subgroup $H$ in $K$, respectively, where $H,K$ are linear algebraic groups over $F$ or $\C$. Note that $N(M,G)=N(A_M,G)$ and $Z(A_M,G)=M$; put $W(M,G):=N(M,G)/M$; likewise for $(M^*,G^*)$, $(\hat M^*,\hat G^*)$, etc.
    \item $N_{\psi_M}(\hat M^*,\hat G^*):=N(A_{\hat M^*},S_{\psi})=\{u\in N(\hat M^*,\hat G^*):\Ad(u)\circ \psi_M = \psi_M\}$. The image of $N_{\psi_M}(\hat M^*,\hat G^*)$ in $W(\hat M^*,\hat G^*)$ is denoted by $W_{\psi_M}(\hat M^*,\hat G^*)$, abbreviated as $W_{\psi_M}$.
    \item $N_{\psi_M}(\hat M^*,\hat G^*)^\natural:=N(A_{\hat M^*},S_\psi)/Z(A_{\hat M^*},S_\psi^{\tu{rad}})$.
\end{itemize}

There is yet another simplification that we will take advantage of, in addition to the two described above. In general the normalized intertwining operator $R_P(u,(\xi,z),\psi_M,\psi_F)$ depends on an element $u \in N_{\psi_M}(\hat M^*,\hat G^*)$, not just its image $w \in W_{\psi_M}$. Write $S_{\psi_M}^{\natural\natural}:=Z(A_{\hat M^*},S_\psi)/Z(A_{\hat M^*},S^{\tu{rad}}_\psi)$. As in \cite[Diagram 2.1.1]{KMSW} there is an exact sequence of groups
\begin{equation} \label{eq:splnpsi}
    1 \to S_{\psi_M}^{\natural\natural} \to N_{\psi_M}(\hat M^*,\hat G^*)^{\natural} \to W_{\psi_M} \to 1.
\end{equation}

The simplification we will use is that this exact sequence splits canonically in our situation. It will therefore be enough to consider $R_P(u,(\xi,z),\psi_M,\psi_F)$ for $u$ in the image of this splitting. We will therefore write $R_P(w,(\xi,z),\psi_M,\psi_F)$ instead, for $w \in W_{\psi_M}$, which is understood to mean the operator obtained for $u$ being the image of $w$ under the splitting.

To explain the splitting, write $\psi_M=\bigoplus_i\psi_i$ according to the blocks of $\hat M^*$. Choose the same homomorphism on blocks with isomorphic parameters.
The Langlands--Shelstad lift $\widehat{\tilde w}$ acts on $\hat M^*$
by the block permutation $w$ by \cite[Lem.~0.4.16]{KMSW}, so $\widehat{\tilde w}\in N_{\psi_M}(\hat M^*,\hat G^*)$.
By the cocycle formula in the proof of \cite[Lem.~2.3.1]{KMSW}, the
multiplicativity defect $\widehat{\tilde w}_1  \widehat{\tilde w}_2 \cdot \widehat{\widetilde{w_1w_2}}^{-1}$ equals, up to conjugation by a block
permutation, the value at $-1$ of a sum of coroots. The roots in
this sum are unions of complete sets of roots between pairs of blocks,
so the sum is a cocharacter of $(Z_{\hat M^*}\cap\SL_N)^\circ$.
This torus lies in $Z(A_{\hat M^*},S_\psi^{\mathrm{rad}})$, so the defect lies there as well. It follows that $w\mapsto \widehat{\tilde w}$ splits \eqref{eq:splnpsi}.

According to \cite[Lem.~0.4.18]{KMSW} we may modify $(\xi,z)$ within its equivalence class so that it satisfies the following.

\begin{itemize}
    \item $\xi(P^*)=P$ and $\xi(M^*)=M$. Then $z$ takes values in $M^*$, $(\xi,z) : M^* \to M$ is an extended pure inner twist, and $\xi$ induces an $F$-isomorphism $W(M^*,G^*) \to W(M,G)$.
    \item $z$ commutes with the Langlands--Shelstad lift of every element of $W(M^*,G^*)(F)$.
\end{itemize}

The isomorphisms above do not depend on the choice of $\xi$ within its equivalence class, subject to the assumption $\xi(M^*)=M$ and $\xi(P^*)=P$, since these assumptions reduce the ambiguity of $\xi$ to an inner automorphism of $M^*(\bar F)$.

Furthermore, by duality of $G^*$ and $\hat G^*$ we have the $\Gamma$-equivariant isomorphism $W(M^*,G^*) \cong W(\hat M^*,\hat G^*)$. Finally, $W_{\psi_M} \subset W(\hat M^*,\hat G^*)^\Gamma$.
Let $w\in W_{\psi_M}$. We transport $w$ to an element of $W(M^*,G^*)(F)$ via the above isomorphisms, which we still denote by $w$; the latter has the Langlands--Shelstad lift $\tilde w\in N(M^*,G^*)(F)$. Let $\breve w := \xi(\tilde w) \in N(M,G)$. It is an $F$-point, since $z$ commutes with $\tilde w$. 
All the linear operators considered below respect the $G$-actions; we will not repeat this point. 

Define a meromorphic family of unnormalized intertwining operators
\[
J_P(\breve{w},\pi_{M,\lambda})\colon  \mathcal{I}_P(\pi_{M,\lambda}) \to \mathcal{I}_{P}(\breve{w}\pi_{M,\lambda})
\]
by meromorphic continuation in the variable $\lambda$ of the following formula:
\[
f\mapsto \left( g\mapsto \int_{(\breve w N \breve{w}^{-1}\cap N)\backslash N} f(\breve{w}^{-1} u  g) du\right),
\]
where we use the measures of \cite[\S2.2]{KMSW}. In the matrix model of \S\ref{ss:notation}, the measures on opposite root spaces are dual under the reduced-trace pairing.
Write $\hat{\mathfrak{n}}:=\Lie \hat{N}^*$. Let $\rho_{w}$ denote the adjoint representation of $^L M^*$ on 
$
 \widehat{\tilde w}^{-1}\hat{\mathfrak{n}}\widehat{\tilde w}/(\widehat{\tilde w}^{-1}\hat{\mathfrak{n}}\widehat{\tilde w}\cap \hat{\mathfrak{n}})
$.
Then we define the normalized operator 
\[
R_P(\breve{w},\psi_{M,\lambda}, \psi_F):=\frac{J_P(\breve{w},\pi_{M,\lambda})}{r_P(w,\psi_{M,\lambda}, \psi_F)}: \mathcal{I}_P(\pi_{M,\lambda}) \to \mathcal{I}_{P}(\breve{w}\pi_{M,\lambda}),
\]
where the normalizing factor is given by Artin--Deligne $L$ and $\epsilon$-factors:
\[
r_P(w,\psi_{M,\lambda}, \psi_F):=\frac{L(0,\rho^\vee_w\circ \phi_{\psi_{M,\lambda}})}{\epsilon(0,\rho^\vee_w\circ \phi_{\psi_{M,\lambda}}, \psi_F)L(1,\rho^\vee_w\circ \phi_{\psi_{M,\lambda}})}.
\]
In the last formula we ignored the Keys--Shahidi factor since it always equals $1$ in our situation.

\begin{thm}\label{thm:regularity}
The normalized intertwining operator $R_P(\breve{w},\psi_{M,\lambda}, \psi_F)$ is regular and invertible at $\lambda=0$.
\end{thm}
\begin{proof}
See Appendix~\ref{akms-reg:appendix}.
\end{proof}

Notice that $\breve{w}\pi_M \cong \pi_M$ since $w \in W_{\psi_M}$.
We recall from \cite[\S2.4.3]{KMSW} how to define the linear operator 
\[
\pi_M(\breve{w},\psi_M)\colon  \breve{w}\pi_M\isom \pi_M,
\]
where $u^\natural$ is the class of the Langlands--Shelstad lift $\widehat{\tilde w}$. Let $\theta=\Ad(\breve w)|_M$ and let
$\hat\theta^*=\Ad(\widehat{\widetilde{w^{-1}}})|_{\hat M^*}$ be its
dual pinned automorphism. By the calculation in
\S\ref{akms-reg:sec:multiplicativity}, we may use
$\widehat{\widetilde{w^{-1}}}^{-1}$ as the representative of the
selected class $u^\natural$.
Then $(M,\theta)$ admits a twisted endoscopic datum $(M',s_{\psi_M})$, 
where $\hat M'$ is the fixed subgroup of $\hat M^*$ by $\Ad(s_{\psi_M})\circ  \hat \theta^*$, with respect to the $L$-morphism $^L M' \hookrightarrow {}^L M^*$ extending $\hat M'\hookrightarrow \hat M^*$ trivially on the Weil group. Then we may assume that $\psi_M\colon \cL_F\times \SL_2(\C)\to {}^L M^*$ factors through $\psi'_M\colon \cL_F\times \SL_2(\C)\to {}^L M'$ possibly after conjugation; write $\pi'_M = \pi_{\psi_M'}\in \Irr_u(M'(F))$ for 
the representation associated to $\psi_M'$.
By \cite[Lem.~2.4.1]{KMSW}, we have a unique isomorphism
$\pi_M(\breve{w},\psi_M)$ such that 
\begin{equation}\label{eq:check-w-characterization}
   \tr \pi'_M(f')=e(M^\theta)\cdot \tr\big(\pi_M(\breve{w},\psi_M)\circ \pi_M(f)\big), 
\end{equation}
whenever $f\in \mathcal H(M)$ transfers to $f'\in \mathcal H(M')$ for the twisted endoscopic datum above. 
Note that $\pi_M(\breve{w},\psi_M)$ depends not only on $\pi_M$ but also on $\psi_M$ since $\pi'_M$ depends on $\psi_M$. (Recall that $\pi_M$ does not determine $\psi_M$ in general, cf.~Remark \ref{rem:not-inj-surj}.)
See Example \ref{ex:weyl-action} below for a special case of \eqref{eq:check-w-characterization} made concrete. 

Finally we define the composite operator
\begin{equation}\label{eq:RP(w)}
    R_P(w, (\xi,z), \psi_{M}, \psi_F)\colon   
\mathcal{I}_P(\pi_{M}) \xrightarrow{R_P(\breve{w},\psi_{M}, \psi_F)} \mathcal{I}_{P}(\breve{w}\pi_{M})\xrightarrow{\mathcal{I}_P(\pi_M(\breve{w},\psi_M))}  \mathcal{I}_P(\pi_{M}).
\end{equation}
Since $\mathcal{I}_P(\pi_{M})$ is irreducible, $R_P(w, (\xi,z), \psi_{M}, \psi_F)$ is a scalar operator. By slight abuse of notation we will write $R_P(w, (\xi,z), \psi_{M}, \psi_F)=c$ for a number $c\in \C$.

It is not hard to see that the operator \eqref{eq:RP(w)} coincides with the one in \cite[(2.5.1)]{KMSW} (via the canonical splitting of \eqref{eq:splnpsi}), taking the following three points into account. First, we compressed the composite of the first two arrows in their formula into a single map. Second, the Keys--Shahidi constant $\lambda(w,\psi_F)$ equals $1$ since $G^*$ is split. Third, our definition already incorporates the cancellation of a certain root number (denoted by $\epsilon_P(w,\psi,\psi_F)$ in \cite[(2.3.1)]{KMSW}), which is multiplied in $l_P(w,\xi,\psi,\psi_F)$ and then divided in $R_{w^{-1}P|P}(\xi,\psi)$ in the notation of \cite[(2.5.1)]{KMSW}.

\begin{lem}\label{lem:independent}
The definition of $R_P(w,(\xi,z),\psi_{M}, \psi_F)$ depends on $(G,\xi,z)$ only through the isomorphism class of $(G,\xi)$ as an inner twist.
Moreover, $R_P(w,(\xi,z),\psi_{M}, \psi_F)$ is independent of the choice of a nontrivial unitary additive character $\psi_F \colon F \rightarrow \C^\times$.
\end{lem}

Henceforth we denote $R_P(w,(\xi,z),\psi_{M}, \psi_F)$ by $R_P(w, \psi_M)$ for simplicity.

\begin{proof}
The first assertion is \cite[Lem.~2.6.4]{KMSW} (more precisely, the equality $\langle c,u^\natural\rangle=1$  in its proof).
For the dependence on $\psi_F$, 
another nontrivial unitary additive character $\psi_F'$ is written as 
$\psi_F'(x) = \psi_F(ax)$ for some $a \in F^\times$.
If we replace $\psi_F$ with $\psi_F'$, 
then 
\begin{itemize}
    \item the invariant measure $du$ is replaced with 
    $|a|^{\frac{1}{2}l}du$
    where 
    \[l = \dim((\widehat{\tilde w}^{-1}\hat{\mathfrak{n}}\widehat{\tilde w}\cap \hat{\mathfrak{n}})\backslash \widehat{\tilde w}^{-1}\hat{\mathfrak{n}}\widehat{\tilde w}),\]
    \item the representative $\breve{w}$ is not changed, 
    \item the normalizing factor $r_P(w,\psi_{M,\lambda},\psi_F)$
    is replaced with 
    \[
    \det(\rho_w^\vee \circ \phi_{\psi_{M,\lambda}})(a)^{-1} |a|^{\frac{1}{2}l}r_P(w,\psi_{M,\lambda},\psi_F),
    \]
    \item the linear operator $\pi_M(\breve{w}, \psi_M)$ is unchanged: it depends on the transfer $f \mapsto f'$, which is normalized by Whittaker data, but in our case there is a unique Whittaker datum. 
\end{itemize}
Order the blocks by $P$, and put $A_i=\phi_{\psi_i}$, $n_i=\dim A_i$,
and $\delta_i=\det A_i$. The determinant at $\lambda=0$ is
\[
 \det(\rho_w^\vee\circ\phi_{\psi_M})
 =\prod_{i<j,\,w(i)>w(j)}\delta_i^{n_j}\delta_j^{-n_i}
 =\frac{\prod_{i<j}\delta_i^{n_j}}{\prod_{w(i)<w(j)}\delta_i^{n_j}}
 =1.
\]
The last equality follows by reindexing, since $A_{w(i)}\simeq A_i$.
Thus $R_P(w,\psi_M)$ is independent of $\psi_F$.
\end{proof}

\begin{lem}\label{lem:mult}
If $w=w'w''$ in $W_{\psi_M}$ then $R_P(w,\psi_{M})=R_P(w',\psi_{M})R_P(w'',\psi_{M})$. 
\end{lem}

\begin{proof}
This follows from \cite[Lem.~2.5.3]{KMSW} but the latter relies on \cite[Lem.~2.2.4]{KMSW} whose proof was postponed for non-tempered parameters.
We complete the proof in Corollary \ref{akms-reg:cor:compound-multiplicativity} below.
\end{proof}

\begin{ex}\label{ex:weyl-action}
Let us explicate \eqref{eq:check-w-characterization} when $G=\GL_m(D)$, $M=\GL_{m/2}(D)\times \GL_{m/2}(D)$, $\psi_M=(\psi^\circ_M,\psi^\circ_M)$, and $\breve w\in N(M,G)$ denotes the Langlands--Shelstad representative of 
the unique nontrivial element in $W(M,G)$ 
permuting the two $\GL_{m/2}(D)$-factors. Assume that $\psi^\circ_M$ is $D$-compatible and let $\sigma\in \Irr_u(\GL_{m/2}(D))$ and $\sigma^*\in \Irr_u(\GL_{dm/2}(F))$ denote the representations corresponding to $\psi^\circ_M$. 
Then $\pi_M=\sigma\otimes \sigma$, $\pi'_M=\sigma^*$. We identify $\theta\colon (x_1,x_2)\mapsto (x_2,x_1)$ in $\Aut_F(M)$; likewise $\hat \theta^*$ is the swapping action on $\hat M^*$. The twisted endoscopic datum $(M',s_{\psi_M})$ for $(M,\theta)$ is isomorphic to $(\GL_{dm/2},1)$ since $g s_{\psi_M}\hat\theta^*(g)^{-1}=1$ for $g=(s_{\psi_M^\circ},1)$.

For a decomposable function $f=f_1\otimes f_2\in \mathcal H(M)$, its twisted transfer $f'\in \mathcal H(\GL_{dm/2}(F))$ can be given by an untwisted transfer (between inner twists) of the convolution $f_1*f_2\in \mathcal H(\GL_{m/2}(D))$. 
Since $(x_1,x_2)\in M(F)$ acts on $\breve w \pi_M=\breve w(\sigma\otimes \sigma)$ by $\sigma(x_2)\otimes\sigma(x_1)$, we can write $\pi_M(\breve{w},\psi_M)= c A$ for a unique constant $c\in \C$, where $A(v_1\otimes v_2)=v_2\otimes v_1$, $v_1,v_2\in \sigma$. Now \eqref{eq:check-w-characterization} becomes
\[
\tr \sigma^*(f') = e(\GL_{m/2}(D))\cdot c \cdot \tr (A\circ (\sigma(f_1)\otimes \sigma(f_2))).
\]
Thanks to \cite[\S1.1.3]{KMSW}, this simplifies as 
\[
\tr \sigma^*(f') = c \cdot e(\GL_{m/2}(D)) \cdot \tr \sigma(f_1*f_2).
\]
Comparing with Proposition \ref{prop:local-JL}, we deduce that $c=a_{\psi_M^\circ}\in \{\pm1\}$.
In conclusion,
\begin{equation}\label{eq:w-action-vs-permutation}
\pi_M(\breve{w},\psi_M)\colon \sigma\otimes \sigma\to \sigma\otimes \sigma,\qquad v_1\otimes v_2\mapsto a_{\psi_M^\circ}\cdot (v_2\otimes v_1).
\end{equation}
\end{ex}

\subsection{The statement of LIR}

We maintain the notation from the preceding subsection. The rest of the paper is devoted to proving the following main theorem.

\begin{thm}\label{thm:LIR}
Given $G$, $G^*$, $P$, $M$, $\psi_M$ as above, we have
\[
\tag{LIR}
R_P(w,\psi_M)=1
\]
for all $w \in W_{\psi_M}$. 
\end{thm}

\begin{lem}\label{lem:easy-case}
If either $w=1$ or $G=\GL_N$, 
then $R_P(w,\psi_M)=1$.
\end{lem}

\begin{proof}
The case $w=1$ follows from invertibility and multiplicativity of $w\mapsto R_P(w,\psi_M)$ (Theorem~\ref{thm:regularity}, Lemma~\ref{lem:mult}). When $G=\GL_N$ this is \cite[Thm.~3.5.1]{AGIKMS}. Namely, if we use identical standard
modules and Whittaker functionals on blocks with equal parameters, the permutation of the standard modules preserves the product Whittaker functional and descends to the permutation of factors on the Langlands quotients; thus the isomorphism in that theorem agrees with \eqref{akms-reg:eq:signed-permutation} for $\GL_N$ (in which case $\kappa(w)$ there equals $1$).

\end{proof}

Let $\zeta\colon W_F\ra Z(\hat G^*)$ denote a continuous character with bounded image. It gives rise to a continuous unitary character
$\chi_\zeta\colon G(F)\to \C^\times$, thus also a continuous unitary character $M(F)\to \C^\times$ upon restriction; the latter will also be denoted by $\chi_\zeta$. Notice that $\psi_M\otimes\zeta\in \Psi(M^*)$ and that $W_{\psi_M}=W_{\psi_M\otimes \zeta}$. Since $\pi_M = \pi_{\psi_M}$,
the compatibility of the Jacquet--Langlands correspondence with character twists tells us that 
$\pi_{\psi_M \otimes \zeta} = \pi_M \otimes \chi_\zeta$.

\begin{lem}\label{lem:character-twist}
Let $\zeta$ be as above. Then $R_P(w,\psi_M)=R_P(w,\psi_M\otimes \zeta)$. In particular $R_P(w,\psi_M)=1$ if and only if $R_P(w,\psi_M\otimes \zeta)=1$.
\end{lem}

\begin{proof}
It follows from the definition that each arrow in \eqref{eq:RP(w)} is compatible with the character twist by $\chi_\zeta$; note that $\breve{w}(\pi_M\otimes\chi_\zeta)=\breve{w}(\pi_M)\otimes\chi_\zeta|_{M(F)}$ since $\chi_\zeta$ is a character of $G(F)$. To verify that $\chi_\zeta$ does not contribute a scalar factor to the intertwining integral, it suffices to note that $\chi_\zeta(\breve{w})=1$ since $\breve{w}$ has reduced norm $1$; the latter holds because the Langlands--Shelstad lift $\tilde{w}$ in $G^*$ has $\det(\tilde{w})=1$.
\end{proof}

\subsection{Initial reductions}

\begin{lem}\label{lem:reduction-to-maximal-Levi}
Fix $F$ and $D$. To prove Theorem~\ref{thm:LIR} it suffices to take
$G=\GL_{2h}(D)$ for all $h\geq 1$, with $P$ the standard parabolic
with Levi $M=\GL_h(D)\times\GL_h(D)$, $\psi_M=(\eta,\eta)$ for a
$D$-compatible $\eta\in\Psi_2(\GL_{dh})$, and $w$ exchanging the two
factors.
\end{lem}

\begin{proof}
Fix $w\in W_{\psi_M}$. If $\psi_M\notin\Psi_2(M^*)$, Lemma~\ref{akms-ind:lem:refinement}
identifies $R_P(w,\psi_M)$ with $R_{P_0}(w_0,\psi_0)$ for
$\psi_0\in\Psi_2(M_0^*)$. Thus we may assume $\psi_M\in\Psi_2(M^*)$.
The group $W_{\psi_M}$ is generated by transpositions of blocks with
equal parameters. By Lemmas~\ref{lem:mult} and \ref{lem:easy-case}, it
suffices to take $w=(i\,j)$ with $\psi_i=\psi_j=\eta$.

Put $h=m_i=m_j$ and $H=\GL_{2h}(D)$, and let $Q$ be its standard
parabolic of type $(h,h)$. Let $S$ be the standard parabolic of
$G^+=\GL_{m+2h}(D)$ with Levi $G\times H$. Using the block order of
$P$ in the first $m$ coordinates, put $P^+=(P\times Q)N_S$ and
$\psi^+=(\psi_M,\eta,\eta)$. Extend $w$ to fix the two new blocks
and let $s$ exchange them. They are conjugate in $W_{\psi^+}$, so
their scalars in $G^+$ agree by Lemma~\ref{lem:mult}.

We claim that
\[
 R_P(w,\psi_M)=R_Q^H(s,(\eta,\eta)).
\]
To verify the claim, compare the operators in $G^+$ with those in
$G$ and $H$ by normalized induction in stages through $S$.
The Langlands--Shelstad lifts of $w$ and $s$ are the lifts in
$G_{\mathrm{der}}$ and $H_{\mathrm{der}}$ embedded in $G^+$. They
normalize $N_S$ and lie in $\ker\delta_S$. The $N_S$ factor of
$N_{P^+}=(N_P\times N_Q)N_S$ therefore cancels from the quotients
in the intertwining integrals of \S\ref{ss:normalized-intertwining-operators},
leaving the quotients and measures for $P$ in $G$ and $Q$ in $H$,
respectively. The same comparison of the dual Lie algebra quotients
identifies the normalizing factors. This identifies the first arrows in
\eqref{eq:RP(w)} by meromorphic continuation and Theorem~\ref{thm:regularity}.
Formula~\eqref{akms-reg:eq:signed-permutation} identifies the second arrows,
since each transposition fixes the blocks in the other factor.
Thus the two equal scalars in $G^+$ are the two sides of the claimed equality.
\end{proof}

We treat the cases in Lemma~\ref{lem:reduction-to-maximal-Levi}
in \S\ref{s:proof-tempered} and \S\ref{s:proof-non-tempered},
according as $\eta$ is tempered or non-tempered. By
Lemma~\ref{lem:character-twist}, we may twist $\eta$ by a unitary character.

\subsection{Extension of LIR beyond $A$-parameters}\label{ss:extension-beyond-A-parameters}

Although LIR is naturally formulated for $A$-parameters, we extend it a little further since it is not guaranteed that we stay within the set of $A$-parameters due to a possible failure of the generalized Ramanujan conjecture for general linear groups. 
In fact, it is enough to extend Lemma \ref{lem:easy-case}, which asserts LIR for $G=\GL_N$.

So assume $G=G^*$, $M=M^*$, and consider a more general parameter $\psi_M\in \Psi^+_u(M^*)$. The earlier construction of the operator $R_P(w, (\xi,z), \psi_{M}, \psi_F)$ goes through without changes, except that it is a priori unknown to be regular or unitary. Since the parabolic induction of $\pi_{\psi_M}\in \Irr_u(M)$ is irreducible, it is a scalar operator if regular.
Lemma \ref{lem:easy-case} can be extended as follows.

\begin{lem}\label{lem:LIR-GL-beyond-A-parameters}
In the above notation $R_P(w, (\xi,z), \psi_{M}, \psi_F)=1$.
\end{lem}

\begin{proof}
Write $\psi_M=(\psi_k)_{k=1}^r$ and, as in \eqref{eq:psi-decomposition},
$\psi_k=\bigoplus_\alpha\eta_{k,\alpha}|\cdot|^{\epsilon_{k,\alpha}}$
with $\eta_{k,\alpha}$ irreducible and $\epsilon_{k,\alpha}\in \R$ with $|\epsilon_{k,\alpha}|<\tfrac12$.
Choose the same ordered decomposition for $\psi_k$'s in the same $w$-orbit. Replacing
$\epsilon_{k,\alpha}$ by $t\epsilon_{k,\alpha}$, we define a deformed parameter $\psi_{M,t}$, which is still fixed
by $w$. Let $\pi_{M,t}$ be the normalized induction on $M$ of the
corresponding twists of Speh representations. To form the family in \eqref{eq:RP(w)},
use $\cI_P(T_w)$ as the second arrow for all $t$, where $T_w$ is the $w$-permutation as defined in \S\ref{akms-reg:sec:multiplicativity}. When $\psi_{M,t}\in\Psi_u^+(M^*)$, the map $T_w$ agrees with $\pi_{M,t}(\breve w,\psi_{M,t})$ by the calculation giving  \eqref{akms-reg:eq:signed-permutation} for split $G$.

With the additional twists $\lambda=(\lambda_k)_{k=1}^r$ of
\S\ref{ss:normalized-intertwining-operators}, the exponents are
$s_{k,\alpha}=t\epsilon_{k,\alpha}+\lambda_k$. As in the proof of Lemma~\ref{akms-ind:lem:refinement}, induction in stages identifies the intertwining integral and its normalizing factor with those for induction from twists of the Speh representations.
By \cite[Prop.~I.10]{MW89}, the normalized intertwining operator is holomorphic where $|\operatorname{Re}(s_{k,\alpha})|<\tfrac12$ for all $k,\alpha$.
For $\lambda=0$, the operator in \eqref{eq:RP(w)} is the identity on $i\R$ by Lemma~\ref{lem:easy-case}, hence at $t=1$ by analytic continuation. Specializing the holomorphic family at $(t,\lambda)=(1,0)$, we obtain  $R_P(w,\psi_M)=1$ as desired.
\end{proof}

\section{Globalization}
\label{s:globalization}

\subsection{The global intertwining relation}
 
Let $\dot F$ be a number field. We fix a nontrivial unitary additive character $\psi_{\dot F}:\dot F\backslash \A_{\dot F}\to \C^\times$, which gives rise to local characters $\psi_{\dot F_v}:\dot F_v\to \C^\times$ and normalized Haar measures on $\dot F_v$ at all places $v$.
With the notation of \S\ref{ss:globalJL}, let $\Psi(\GL_{N,\dot F})$
be the set of formal sums $\boxplus_i\dot\psi_i$ with
$\dot\psi_i\in\Psi_2(\GL_{N_i,\dot F})$ and $\sum_iN_i=N$.
Extend \eqref{eq:localize-A-parameters} by localizing each summand and
taking their direct sum. For a split Levi subgroup, $\Psi_2$ and
localization are defined factor by factor.

Let $\dot D$ be a central division algebra over a number field $\dot F$. Set $d:=[\dot D:\dot F]^{1/2}$ and $\dot D_v:=\dot D\otimes_{\dot F}\dot F_v$. Let $\tu{Ram}(\dot D)$ denote the set of places of $\dot F$ where $\dot D$ is ramified.
Consider proper Levi and parabolic subgroups 
\[
\dot M=\textstyle\prod_{i=1}^r  \GL_{m_i}(\dot D)\subset \dot P\subset \dot G=\GL_m(\dot D)
\]
over $\dot F$, where the parabolic subgroup $\dot P$ is block upper triangular with $\dot M$ as a Levi component. We have the corresponding Levi and parabolic subgroups 
\[
\dot M^*=\textstyle\prod_{i=1}^r  \GL_{d m_i} \subset \dot P^*\subset  \dot G^*=\GL_N,\qquad  N=dm.
\]
Let $\dot \psi_M\in \Psi_2(\dot M^*)$, so $\dot\psi_{M,v}\in \Psi^+_u(\dot M^*_v)$. Write $\dot \pi_{M,v}^*$ for the unique member of $\Pi_{\dot \psi_{M,v}}(\dot M_v^*)$, and put $\dot \pi^*_M:=\otimes'_v \dot \pi^*_{M,v}$.
For the inner form $\dot M$, at every place $v$ of $\dot F$ where $\dot D_v$ is not split (equivalently where $\dot M_v$ is not split), suppose that 
\begin{equation}\label{eq:packet-nonempty-at-v}
    \dot \psi_{M,v}\in \Psi(\dot M^*_v)_{\dot D_v\tu{-comp}},\qquad \forall v\in \tu{Ram}(\dot D). 
\end{equation}
At these places let $\dot\pi_{M,v}$ be the unique member of
$\Pi_{\dot\psi_{M,v}}(\dot M_v)$. At split places put
$\dot\pi_{M,v}=\dot\pi^*_{M,v}$ under the identification of the
two groups. Put $\dot\pi_M=\otimes'_v\dot\pi_{M,v}$.
By Proposition~\ref{prop:global-JL}, applied to each factor of $\dot M$,
$\dot\pi_M$ occurs in the discrete spectrum with multiplicity one.

Define $W_{\dot\psi_M}=W_{\dot\psi_M}(\dot M^*,\dot G^*)$
as in the local case and let $W^{\tu{reg}}_{\dot\psi_M}$ be the
subset of elements $w$ with
$\det(w-1|\fka_{\dot M^*}/\fka_{\dot G^*})\ne0$.
Define the global operator
\[
R_{\dot P}(w,\dot \psi_M):=\prod_v R_{\dot P_v}(w,\dot \psi_{M,v}), \qquad w\in W_{\dot\psi_M}, 
\]
on the normalized parabolic induction of $\dot\pi_M$ to $\dot G$.
The factors at split places are the identity by
Lemma~\ref{lem:LIR-GL-beyond-A-parameters}, so the product is finite.
Define $R_{\dot P^*}(w,\dot\psi_M)$ in the same way for $\dot\pi_M^*$.
The same lemma gives
\begin{equation}\label{eq:R(w)=1-GL(N)}
    R_{\dot P^*}(w,\dot \psi_M)=1.
\end{equation}

To verify the analogue for non-split inner forms, we start by considering the (un-normalized) global intertwining operators $M_{\dot P^*}(w,\dot\pi^*_M)$ and $M_{\dot P}(w,\dot\pi_M)$ on $\cI_{\dot P^*}(\dot \pi^*_M)$ and $\cI_{\dot P}(\dot \pi_M)$, respectively, in the trace formula context, as introduced in \cite[(3.1.6)]{Arthur} (restricted to induced representations of $\dot \pi^*_M$ and $\dot \pi_M$, cf.~\cite[p.182]{Arthur}).
\begin{lem}\label{lem:Global-Intertwining-Relation}
Let $\dot F, \dot D, \dot G, \dot P, \dot M, \dot\psi_M$ be as above. Suppose that $W^{\tu{reg}}_{\dot\psi_M}$ is a singleton, and write $W^{\tu{reg}}_{\dot\psi_M}=\{w\}$. Then the following equality of scalars holds.
\[
M_{\dot P^*}(w,\dot\pi^*_M)=M_{\dot P}(w,\dot\pi_M).
\]
\end{lem}

\begin{proof}
In this proof, we freely use the standard notation for Arthur's invariant trace formula as in \cite[\S3.1, \S3.3]{Arthur}. In particular, we use the global Hecke algebra $\mathcal H(\dot G)$ for $\dot G$ which is translation-invariant by the central subgroup $A_{\dot G,\infty}^+$ of $\dot G(\dot F \otimes_{\Q} \R)$ in the notation there; $\mathcal H(\dot G^*)$ is defined likewise. 
By surjectivity of \eqref{eq:B-to-Bad}, we promote $\dot G$ to a global extended pure inner twist of $\dot G^*$; we fix this structure.
Write $\dot\psi\in \Psi(\dot G^*)$ for the image of $\dot \psi_M$.

For each $t\ge 0$, we have the equality of invariant trace formulas for $\dot G$ and $\dot G^*$ from \cite[(25.28)]{Art05}, which is one of the main results in \cite{AC89}:
\[
I^{\dot G}_{\disc,t}(\dot f)=I^{\dot G^*}_{\disc,t}(\dot f^*), \qquad \dot f\in \mathcal H(\dot G),\quad \dot f^*\in \mathcal H(\dot G^*),
\]
where $\dot f=\prod_v \dot f_v$, $\dot f^*=\prod_v \dot f_v^*$, and $\dot f^*_v$ is a transfer of $\dot f_v$ at every place $v$. 

The global parameter $\dot \psi$ determines $t=t(\dot\psi)\ge 0$ as well as a family of Satake parameters for $\dot G^*$. Write $c_{\dot G^*}(\dot\psi)$ for the family of Satake parameters thus obtained. 
This allows us to define $I^{\dot G}_{\disc,\dot \psi}$ and $I^{\dot G^*}_{\disc,\dot \psi}$ by separating the representations whose families of Satake parameters match $c_{\dot G^*}(\dot\psi)$. Then the above equality is refined as
\begin{equation}\label{eq:I-disc-psi=I-disc-psi-*}
I^{\dot G}_{\disc,\dot \psi}(\dot f)=I^{\dot G^*}_{\disc,\dot \psi}(\dot f^*).
\end{equation}
See \cite[\S3.3]{Arthur} for details. Following p.145 and pp.188--189 in \emph{op.~cit.}, we can express $I^{\dot G^*}_{\disc,\dot \psi}(\dot f^*)$
in the following form:
\begin{equation}\label{eq:I-disc-psi}
    \sum_{\dot L^*} |W(\dot L^*)|^{-1} \sum_{w\in W^{\tu{reg}}(\dot L^*)} \frac{\tr\big(M_{\dot Q^*,\dot\psi}(w)\cI_{\dot Q^*,\dot\psi}(\dot f^*)\big)}{\big|\det(w-1| \fka_{\dot L^*}/\fka_{\dot G^*})\big|},
\end{equation}
where 
\begin{itemize}
    \item $\dot L^*$ runs over a set of representatives for $\dot G^*(\dot F)$-conjugacy classes of Levi subgroups of $\dot G^*$; we choose $\dot L^*$ to be a block diagonal subgroup,
    \item $W(\dot L^*)$ is the Weyl group of $\dot L^*$ in $\dot G^*$, and $W^{\tu{reg}}(\dot L^*)$ is the subset of $w$ such that the denominator above is nonzero,
    \item $\dot Q^*$ is the block upper triangular parabolic subgroup of $\dot G^*$ with Levi factor $\dot L^*$,
    \item $\cI_{\dot Q^*,\dot\psi}$ is the subspace in the normalized parabolic induction from the $L^2$-discrete spectrum of $\dot L^*$ whose family of Satake parameters is given by $c_{\dot G^*}(\dot\psi)$,
    \item $M_{\dot Q^*,\dot\psi}(w)$ is the (un-normalized) global intertwining operator on $\cI_{\dot Q^*,\dot\psi}$ given by $w$. (If no irreducible subrepresentation is stabilized by $w$ then the trace of the corresponding operator vanishes.)
\end{itemize}
A theorem of Jacquet--Shalika implies that $c_{\dot G^*}(\dot\psi)$ contributes to $\cI_{\dot Q^*,\dot\psi}$ only if $\dot L^*$ is conjugate to $\dot M^*$, thus we may concentrate on $\dot L^*=\dot M^*$, and only if the inducing representation is conjugate to $\dot \pi^*_M$ which was introduced above \eqref{eq:packet-nonempty-at-v}. The regular element $w$ in the hypothesis permutes the blocks in a
single cycle. Since it fixes
$\dot\psi_M$, the inducing factors are all equal and the Weyl orbit of
$\dot\pi_M^*$ is a singleton. Hence \eqref{eq:I-disc-psi} reduces to 
\[
|W(\dot M^*)|^{-1}\big|\det(w-1| \fka_{\dot M^*}/\fka_{\dot G^*})\big|^{-1} \tr\big(M_{\dot P^*}(w,\dot\pi_M^*)\cI_{\dot P^*}(\dot \pi_M^*,\dot f^*)\big),
\]
where $w$ is the unique element of $W_{\dot \psi_M}^{\tu{reg}}$, and $\cI_{\dot P^*}(\dot \pi_M^*,\dot f^*)$ denotes the action of $\dot f^*$ defined on the induced representation from $\dot\pi_M^*$ on $\dot M^*(\A_{\dot F})$. 

We have the expansion for $I^{\dot G}_{\disc,\dot\psi}(\dot f)$ analogous to \eqref{eq:I-disc-psi}. By Proposition~\ref{prop:global-JL}, the global Jacquet--Langlands
transfer is injective and matches the Satake parameters at unramified
places. The same argument therefore reduces
$I^{\dot G}_{\disc,\dot\psi}(\dot f)$ to the corresponding term
\[
|W(\dot M)|^{-1}\big|\det(w-1| \fka_{\dot M}/\fka_{\dot G})\big|^{-1} \tr\big(M_{\dot P}(w,\dot\pi_M)\cI_{\dot P}(\dot \pi_M,\dot f)\big),
\]
where $w$ is viewed as an element of $W^{\tu{reg}}(\dot M)$ via the natural identification $W(\dot M)=W(\dot M^*)$. 
Hence \eqref{eq:I-disc-psi=I-disc-psi-*} becomes
\[
\tr\big(M_{\dot P^*}(w,\dot\pi_M^*)\cI_{\dot P^*}(\dot \pi_M^*,\dot f^*)\big)=
\tr\big(M_{\dot P}(w,\dot\pi_M)\cI_{\dot P}(\dot \pi_M,\dot f)\big).
\]
The intertwining operators are scalars since parabolic induction from an irreducible unitary representation is irreducible. By Convention~\ref{convention:transfer-factor}, applied at every place, the induced characters agree provided that $\prod_v \bigl(a_{\dot\psi_v}e(\dot G_v)\bigr)=1$. 
Once this product identity is established, choose $\dot f$ for which $\tr\cI_{\dot P}(\dot\pi_M,\dot f)\neq 0$. Cancelling this common character value in the preceding trace identity proves the lemma.
Since $\prod_v e(\dot G_v)=1$ by \cite{Kot83}, it remains to verify
$
\textstyle\prod_v a_{\dot\psi_v}=1.
$

For this, write $\dot \psi_M=(\dot\psi_{1},...,\dot\psi_{r})$ with $\dot\psi_{i}=\dot\phi_i\boxtimes S_{b_i}\in \Psi_2(\GL_{dm_i})$, so that $\dot\psi=\boxplus_{i=1}^r \dot\psi_i$. Then it suffices to show that $\prod_v a_{\dot\psi_{i,v}}=1$.
Write $s_{\dot\psi_i}\in Z(\hat \GL_{dm_i})=\C^\times$ for the image of $-1\in \SL_2(\C)$ under $S_{b_i}$. 
Let $\dot\pi_i$ denote the discrete automorphic representation of $\GL_{m_i}(\dot D)$ corresponding to $\dot\psi_i$.
Then
\[
\textstyle\prod_v a_{\dot\psi_{i,v}}=\textstyle\prod_v \langle s_{\dot\psi_{i}},\dot\pi_{i,v}\rangle = 1
\]
since $\prod_v \langle \cdot,\dot\pi_{i,v}\rangle$ is trivial on $Z(\hat \GL_{dm_i})$ by Kottwitz's product formula \cite[(0.3.2)--(0.3.3)]{KMSW} (cf.~the penultimate paragraph above Theorem 1.7.1 there); here the pairing is computed with respect to the localization at $v$ of the fixed global extended pure inner twist.
\end{proof}

\begin{cor}\label{cor:global-LIR}
In the setting of the preceding lemma, assume that $r=2$ and $m_1=m_2=1$ so that $\dot M=\dot D^\times\times \dot D^\times$ and $\dot\psi_M$ has the form $(\dot\psi_M^\circ,\dot\psi_M^\circ)$. Further, at each place $v$ of $\dot F$ where $\dot D_v$ is not split, we assume that $\dot\psi_{M,v}\in \Psi_2(\dot M^*_v)_{\dot D_v\tu{-comp}}$. Then
\[
R_{\dot P}(w,\dot \psi_M)=1.
\]
\end{cor}

\begin{proof}
In view of \eqref{eq:R(w)=1-GL(N)}, it suffices to prove that $R_{\dot P^*}(w,\dot \psi_M)=R_{\dot P}(w,\dot \psi_M)$.
This follows from Lemma \ref{lem:Global-Intertwining-Relation} because $R_{\dot P^*}(w,\dot \psi_M)$ and $R_{\dot P}(w,\dot \psi_M)$ are obtained from the two sides of the identity there by dividing by the same global normalizing factor as in \cite[(4.2.10)]{Arthur} (which is regular and nonzero at zero by the Rankin--Selberg functional equation), except that the following sign verification is needed. Namely, since the $w$-action given by $\pi_{\dot M}(\breve w,\dot\psi_{M,v})$ in \eqref{eq:w-action-vs-permutation} differs from the permutation action $v_1\otimes v_2\mapsto v_2\otimes v_1$ by a factor of $a_{\dot\psi^{\circ}_{M,v}}$ at each place $v$ (we have $a_{\dot\psi^{\circ}_{M,v}}=1$ unless $\dot D$ ramifies at $v$), we need to verify that $\prod_v a_{\dot\psi^{\circ}_{M,v}}=1$. This is shown by exactly the same argument as in the last paragraph in the proof of Lemma \ref{lem:Global-Intertwining-Relation}.
\end{proof}

\begin{rem}
The assumption in the corollary is far from optimal but enough for our purpose.
\end{rem}

\subsection{A globalization lemma}

For each local field $F$, define an equivalence relation on $\Irr(\GL_n(F))$ by $\pi_1\sim \pi_2$ if there exists a unitary $\chi \in \Irr(\GL_1(F))$ such that $\pi_1\cong \pi_2\otimes (\chi\circ \det)$. Likewise $\psi_1\sim\psi_2$ in $\Psi(\GL_n(F))$ if $\psi_1\cong \psi_2\otimes\chi$ for some $\chi\in \Irr(\GL_1(F))$, where $\chi$ is viewed as a continuous character $W_F\to Z(\hat{\GL_n})=\C^\times$ via local class field theory. These equivalence relations extend to a finite product of general linear groups in the obvious manner.
We can always find a number field $\dot F$ and a place $v_0$ with $\dot F_{v_0}\cong F$, cf.~\cite[Lem.~6.2.1]{Arthur}.

\begin{lem}\label{lem:globalization-rep}
Let $F$ be a local field of characteristic zero and
$\pi\in\Irr_2(\GL_n(F))$. Fix $\dot F$ and $v_0$ with $\dot F_{v_0}\cong F$ as above. There exists
$\dot\pi\in\cA_{\cusp}(\GL_{n,\dot F})$ with $\dot\pi_{v_0}\sim\pi$.
Given distinct finite places $v_1,\dots,v_s$ of $\dot F$ different from
$v_0$, we may impose either of the following conditions on $\dot\pi$.
\begin{itemize}
    \item[(i)] $\dot\pi_{v_i}$ is a unitary unramified character twist of the Steinberg representation for every $i=1,...,s$, or
    \item[(ii)] $\dot\pi_{v_i}$ is unitary supercuspidal for every $i=1,...,s$.
\end{itemize}
\end{lem}

\begin{proof}
Take $S$ to be the finite set of places of $\dot F$ consisting of $v_0, v_1,...,v_s$ and all infinite places. We can deduce from \cite[Thm.~4.8]{Shi12} the existence of $\dot \pi$ which globalizes $\pi$ and satisfies (i) (resp.~(ii)).\footnote{Several results of \cite{Shi12} including Theorem 4.8 rely on \cite{Sauvageot}, which unfortunately contained a number of inaccuracies. They were addressed in the appendix of \cite{Shi12} except for a delicate but crucial point, which was resolved very recently by Takanashi \cite{TakanashiDensity}. See the latter for details.}
More precisely, choose distinct finite places $w_1,w_2$ of $\dot F$ outside $S$. (In the setting of \emph{loc.~cit.}, we prescribe a twist of Steinberg at $w_1$ and a supercuspidal condition at $w_2$.) The reference works with the unitary dual for the quotient group by a central torus $A_{\infty}=\R^\times_{>0}$ in $\GL_n(\dot F\otimes_{\Q}\R)$. Accordingly, we choose an auxiliary subset $\hat U_\infty$ in the tempered dual of $(\prod_{y|\infty}\GL_n(\dot F_y))/A_{\infty}$ with positive Plancherel measure with compact closure and boundary volume zero, with the following requirement if $v_0$ is Archimedean: the $v_0$-component of $\hat U_\infty$ is contained in the unitary character twists of $\pi$. The main point is that both the equivalence class of $\pi$ and condition (i) (resp.~(ii)) define a subset in the unitary dual of $\GL_n(\dot F_S)$ whose Plancherel measure is positive. This is standard; indeed, the locus of unitary unramified twists of a fixed unitary discrete series representation has positive Plancherel measure. (One can also see this from the explicit formula in \cite[\S5]{AubertPlymen}.)
Thus we apply \cite[Thm.~4.8]{Shi12} by taking
\begin{itemize}
    \item $v_1,v_2$ there to be our $w_1,w_2$,
    \item $\hat f_S$ to be the characteristic function of the following sets in the unitary dual of $\GL_n(\dot F_S)/A_{\infty}$:
    \begin{itemize}
    \item the set $\hat U_\infty$ above at the Archimedean places, 
    \item the set of unitary unramified character twists of $\pi$ at $v_0$ if $v_0$ is finite,
        \item at $v_1,...,v_s$, the subset determined by condition (i), resp.~the subset of unitary unramified character twists of a fixed supercuspidal in case of (ii),
    \end{itemize}
\end{itemize}
to see that the right hand side in the formula of \emph{loc.~cit.} is nonzero. The non-vanishing of the left hand side therein immediately implies the desired existence of $\dot \pi$.
\end{proof}

\begin{lem}\label{lem:globalization-parameter} Consider the following data (if $F$ is Archimedean, take $a=1$ and ignore the corresponding $\SL_2$-factor):
\begin{itemize}
    \item $D$ is a central division algebra over $F$ with $d:=[D:F]^{1/2}>1$,
    \item $\psi=\rho\boxtimes S_a\boxtimes S_b\in \Psi_2(\GL_{n,F})$ is $D$-compatible (Definition \ref{def:D-compatible}).
\end{itemize}
Then these data admit a globalization as follows: 
\begin{itemize}
    \item $\dot F$ is a number field,
    \item $v_0$ is a place of $\dot F$, such that $\dot F_{v_0}\cong F$,
    \item $v_1,...,v_s$ are finite places of $\dot F$ different from $v_0$, where $s\in \Z_{\ge1}$,
    \item $\dot D$ is a central division algebra over $\dot F$ unramified away from $v_0,...,v_s$ such that 
    \begin{itemize}
        \item $\dot D_{v_0}\cong M_{n/d}(D)$, so that $\dot D_{v_0}^\times\cong \GL_{n/d}(D)$,
        \item $\dot D_{v_i}$ is a division algebra for $i=1,...,s$.
    \end{itemize}
    \item $\dot\psi\in \Psi_2(\GL_{n,\dot F})$ such that $\dot\psi_{v_0}\sim \psi$, and
    the following hold for all $i=1,...,s$:
    \begin{itemize}
        \item[(i)] in case $b=1$, we have $\dot \psi_{v_i}=\rho_i\boxtimes S_n \boxtimes S_1\in\Phi_{2}(\GL_{n,\dot F_{v_i}})$ with $\dim\rho_i=1$, and
        \item[(ii)] in case $b>1$, we have $\dot \psi_{v_i}=\rho_i\boxtimes S_1\boxtimes S_b\in\Psi_{2}(\GL_{n,\dot F_{v_i}})$.
    \end{itemize}
\end{itemize}
Notice that $\dot\psi_{v_i}$ in (i) and (ii) above are $\dot D_{v_i}$-compatible.
\end{lem}

\begin{rem}
The $D$-compatible condition on $\psi$ ensures that $n/d\in \Z_{\ge1}$.
\end{rem}

\begin{proof}
Write $\tu{inv}(D)=\tfrac{j}{d}\in \tfrac{1}{d}\Z/\Z$, with $1 \leq j < d$ and $(j,d)=1$. Put $s:=(d-j)n/d\in \Z_{\ge 1}$. Fix a number field $\dot F$ and places $v_0,...,v_s$ as in the first three bullet points in the globalization desiderata.
  There exists a central division algebra $\dot D$ over $\dot F$ such that 
\[
\tu{inv}(\dot D_{v})=\begin{cases}
j/d, & v=v_0,\\
1/n, & v=v_1,...,v_s,\\
0, & v\notin \{v_0,...,v_s\},
\end{cases}
\]
since $\sum_{v}\tu{inv}(\dot D_v)=0\in \tfrac{1}{n}\Z/\Z$. By construction $\dot D$ satisfies the conditions in the lemma.

Denote by $\phi\in \Irr_2(\GL_{n/b}(F))$ the representation whose $L$-parameter is $\rho\boxtimes S_a$. By Lemma \ref{lem:globalization-rep} there exists $\dot\phi\in \cA_{\tu{cusp}}(\GL_{n/b,\dot F})$ such that $\dot \phi_{v_0}\sim \phi$ and $\dot \phi_{v}$ satisfies (i) (resp.~(ii)) of that lemma at each $v=v_1,...,v_s$. 
Then $\dot \psi:=\dot\phi\boxtimes S_b$ satisfies (i) (resp.~(ii)) of the current lemma. 
\end{proof}

\section{Proof of LIR: the tempered case}\label{s:proof-tempered}

We treat the tempered case of Lemma~\ref{lem:reduction-to-maximal-Levi},
writing $\phi_M=(\phi_M^\circ,\phi_M^\circ)$ for $\psi_M=(\eta,\eta)$. Thus $N=2dh$.

\subsection{Reduction via globalization}

\begin{lem}\label{lem:reduction-tempered}
In the tempered case, it suffices to prove LIR for the following data.
\begin{itemize}
\item $F$ is non-Archimedean,
    \item $M=\GL_1(D)\times \GL_1(D)$ (in particular $N=2d$),
    \item $\phi_M=(\triv\boxtimes S_d,\triv\boxtimes S_d)$,
    \item $w\in W_{\phi_M}$ is the unique nontrivial element.
\end{itemize}
\end{lem}

\begin{proof}
By Lemma~\ref{lem:easy-case}, we may assume $d>1$.
Assume LIR for the data in the statement.
We globalize $F$, $D$, $\phi_M^\circ$ to $\dot F$, $\dot D$, $\dot\phi_M^\circ$ by Lemma \ref{lem:globalization-parameter}, case (i). 
Hence 
\[\dot \phi_{M,v_i}=(\rho_i\boxtimes S_{N/2},\,\rho_i\boxtimes S_{N/2}),\qquad \dim \rho_i=1,\qquad i=1,...,s.
\]
Put $\dot G=\GL_2(\dot D)$ and let $\dot P$ be its standard parabolic
with Levi $\dot M=\dot D^\times\times\dot D^\times$. Put
$\dot\phi_M=(\dot\phi_M^\circ,\dot\phi_M^\circ)$.
Since $\dot \phi_{M,v_i}=\rho_i\otimes (\triv\boxtimes S_{N/2},\triv\boxtimes S_{N/2})$, our assumption together with Lemma \ref{lem:character-twist} tells us that LIR is true for $(\dot G_{v},\dot M_{v},\dot\phi_{M,v})$ at $v=v_1,...,v_s$. The same holds at $v\notin \{v_0,v_1,...,v_s\}$ by Lemma \ref{lem:LIR-GL-beyond-A-parameters}. On the other hand, 
\[
R_{\dot P}(w,\dot \phi_M)=\textstyle\prod_v R_{\dot P_v}(w,\dot \phi_{M,v})=1
\]
by Corollary \ref{cor:global-LIR}. We have seen that $R_{\dot P_v}(w,\dot \phi_{M,v})=1$ for $v\neq v_0$, so the same holds at $v=v_0$. Therefore we conclude by Lemma \ref{lem:character-twist}.
\end{proof}

\subsection{An explicit computation}
Since $R_P(w,\phi_M)$ is independent of 
the choice of a nontrivial additive character $\psi_F$ of $F$ (Lemma \ref{lem:independent}), 
we may assume hereafter that $\psi_F$ is trivial on the ring of integers $\cO_F$ but not on the inverse of the maximal ideal $\cP_F^{-1}$.
\par

Let $G = \GL_2(D)$. 
We denote by $P = MN$ the upper triangular parabolic subgroup of $G$ so that $M = \GL_1(D) \times \GL_1(D)$
and $N \cong D$. 

\begin{lem}\label{lem:lift-and-vol}
Notations and assumptions are as above.
\begin{enumerate}
    \item The lift of the unique nontrivial element in $W(M,G)$ is given by 
    \[
    \breve{w} = \begin{pmatrix}
        0 & 1 \\ (-1)^d & 0
    \end{pmatrix}.
    \]
    \item
    The measure on $N \cong D$ is determined by the condition $\vol(\cO_D) = q^{-d(d-1)/2}$.
\end{enumerate}
\end{lem}
\begin{proof}
For (i), use the isomorphism $\xi$ of \S\ref{ss:notation}, as allowed by
Lemma~\ref{lem:independent}.
Let $G^* = \GL_{2d}(F)$ and $M^* = \GL_d(F) \times \GL_d(F)$. 
We denote by $T^*$ the diagonal torus of $G^*$.
The unique lift $w_T$ of $w\in W(M^*,G^*)$ to $W(T^*,G^*)$
that preserves the splitting 
$\tu{spl}_{M^*} = (B^* \cap M^*, T^*, \{X_i\}_{1 \leq i \leq 2d-1, i \not= d})$
is the image of $\begin{psmallmatrix}
    0 & 1_d \\ 1_d & 0
\end{psmallmatrix}$ in $W(T^*,G^*)$. 
Here, $1_n$ denotes the identity matrix in $\GL_n(F)$.
For $1 \leq i \leq 2d-1$, 
if we denote by $s_i$
the simple reflection with respect to the simple root 
$T^* \ni t \mapsto t_i/t_{i+1}$, 
then the reduced decomposition of $w_T$ is given by 
\[
w_T = (s_ds_{d+1} \cdots s_{2d-1})
(s_{d-1}s_d \cdots s_{2d-2}) 
\cdots (s_1s_2 \cdots s_d).
\]
If we write 
\[
\tilde{s}_i = \exp(X_i) \exp(-{}^tX_i) \exp(X_i)
= \left(\begin{array}{c|cc|c}
     1_{i-1}&&&  \cr \hline
     &0&1& \cr
     &-1&0& \cr \hline
     &&&1_{2d-i-1}
\end{array}\right), 
\]
then the Langlands--Shelstad lift $\tilde{w}$ of $w$ is given by 
\[
\tilde{w} = (\tilde{s}_d\tilde{s}_{d+1} \cdots \tilde{s}_{2d-1})
(\tilde{s}_{d-1}\tilde{s}_d \cdots \tilde{s}_{2d-2}) 
\cdots (\tilde{s}_1\tilde{s}_2 \cdots \tilde{s}_d)
= \begin{pmatrix}
    0 & 1_d \\ 
    (-1)^d \cdot 1_d & 0
\end{pmatrix}.
\]
Since the inner twist $\xi \colon G^*_{\overline{F}} \xrightarrow{\sim} G_{\overline{F}}$ 
is induced from an algebraic isomorphism 
$M_d(\overline{F}) \xrightarrow{\sim} D \otimes_F \overline{F}$, 
we conclude that 
\[
\breve{w} = \xi(\tilde{w}) = \begin{pmatrix}
0 & 1 \\ (-1)^d & 0
\end{pmatrix}.
\]
\par

We show (ii).
If we denote the reduced trace of $D$ by $\tr_D \colon D \rightarrow F$, 
then 
\[
\{
x \in D \;|\; \tr_D(xy) \in \cO_F,\; \forall y \in \cO_D
\} 
= \cP_D^{-d+1}, 
\]
where $\cP_D$ is the maximal (two-sided) ideal of $\cO_D$.
Let $\Phi$ be the characteristic function on $\cO_D$. 
Then 
\begin{align*}
    \int_D \hat\Phi(x) dx 
    &= \int_D \int_D \Phi(y) \psi_F(\tr_D(xy)) dydx \\
    &= \vol(\cP_D^{-d+1})\vol(\cO_D) \\
    &= q^{d(d-1)}\vol(\cO_D)^2.
\end{align*}
Since we use the self-dual measure, 
the left-hand side is equal to $\Phi(0) = 1$. 
Hence $\vol(\cO_D) = q^{-d(d-1)/2}$.
\end{proof}

Now we prove LIR in the case stated in Lemma \ref{lem:reduction-tempered}. In this case, $\Pi_{\phi_M}=\{\pi_M\}$, where $\pi_M=\triv\boxtimes \triv$ is the trivial representation of $M =D^\times \times D^\times$. 
Slightly deviating from the notation of \S\ref{ss:normalized-intertwining-operators}, we consider a family of twists of $\pi_M$:
\[
\pi_{M,\lambda}:=|\cdot|^{\lambda}\boxtimes |\cdot|^{-\lambda},\quad \lambda\in \C.
\]
Set $K = \GL_2(\cO_D)$.
Write $f_\lambda\in \cI_P(\pi_{M,\lambda})$ for the unique vector such that $f_\lambda(k) = 1$ for $k \in K$.
\par

\begin{lem}\label{lem:unnormalized-intertwining}
The value of $J_{P}(\breve{w},\pi_{M,\lambda})(f_\lambda)\in \cI_{P}(\breve{w}\pi_{M,\lambda})$ at $1\in G$ equals 
\[
\vol(\cO_D^\times)\left(
\frac{q^{-d}}{1-q^{-d}} + \frac{1}{1-q^{-2\lambda}}
\right)
\]
as a meromorphic function in $\lambda$.
\end{lem}

\begin{proof}
We denote by $\mathrm{ord}$ the normalized valuation on $D$. 
Then 
\[
\breve{w}^{-1}\begin{pmatrix}
    1 & x \\ 0 & 1
\end{pmatrix}
=
\begin{pmatrix}
    0 & (-1)^d \\ 1 & x
\end{pmatrix}
= \left\{
\begin{aligned}
&\begin{psmallmatrix}
    1 & 0 \\ 0 & 1
\end{psmallmatrix}
\begin{psmallmatrix}
    0 & (-1)^d \\ 1 & x
\end{psmallmatrix}
&&\text{if $\mathrm{ord}(x) \geq 0$}, \\
&\begin{psmallmatrix}
    x^{-1} & (-1)^d \\ 0 & x 
\end{psmallmatrix}
\begin{psmallmatrix}
    (-1)^{d-1} & 0 \\ x^{-1} & 1
\end{psmallmatrix}
&&\text{if $\mathrm{ord}(x) < 0$}.
\end{aligned}
\right. 
\]
The two cases express the matrix on the left as a product in the Iwasawa decomposition $G=PK$.
Since 
\[
\frac{\vol(\cO_D)}{\vol(\cO_D^\times)} = \frac{1}{1-q^{-d}}, 
\quad 
\frac{\vol(\{\mathrm{ord}(x) = -m\})}{\vol(\cO_D^\times)} = q^{dm},
\]
we have
\begin{align*}
    \frac{1}{\vol(\cO_D^\times)}
    J_{P}(\breve{w},\pi_{M,\lambda})(f_\lambda)(1)
    &= 
    \frac{1}{\vol(\cO_D^\times)} 
    \int_D f_\lambda\left(\breve{w}^{-1}
    \begin{pmatrix}
        1 & x \\ 0 & 1
    \end{pmatrix}\right)
    dx
    \\&= \frac{\vol(\cO_D)}{\vol(\cO_D^\times)}
        + \sum_{m = 1}^\infty\frac{1}{\vol(\cO_D^\times)} \int_{\mathrm{ord}(x) = -m} 
        f_\lambda\left(\begin{pmatrix}
    x^{-1} & (-1)^d \\ 0 & x 
\end{pmatrix}\right) dx
    \\&= \frac{1}{1-q^{-d}} + \sum_{m=1}^\infty q^{-(2\lambda+d)m} \cdot q^{dm}
    \\&= \frac{1}{1-q^{-d}} + \frac{q^{-2\lambda}}{1-q^{-2\lambda}}
    = \frac{q^{-d}}{1-q^{-d}} + \frac{1}{1-q^{-2\lambda}}.
\end{align*}
This completes the proof.
\end{proof}

\begin{lem}\label{lem:LIR-tempered-special-case}
LIR is true in the case specified by Lemma \ref{lem:reduction-tempered}.
\end{lem}

\begin{proof}
In the case at hand, both $\breve{w}\pi_M$ and $\pi_M$ are the trivial representation of $M$, so $\pi_M(\breve{w}, \phi_M)\in \End(\C)=\C$. In fact $\pi_M(\breve{w}, \phi_M)=1$ from Example \ref{ex:weyl-action} (we identify $\sigma\otimes \sigma=\C$ in the case at hand) and Lemma \ref{lem:a-psi}. As a result, in order to show that $R_P(w,\phi_M)=1$, it suffices to check that the function $R_P(\breve{w},\phi_M)(f_0)\in \cI_P(\breve{w}\pi_M)=\cI_P(\triv)$ has value $1$ at $1\in G$.
In light of Lemma \ref{lem:unnormalized-intertwining}, this boils down to verifying that
\[
r_P(w,\phi_{M,\lambda})^{-1}\cdot\frac{\vol(\cO_D^\times)}{1-q^{-2\lambda}} \bigg\rvert_{\lambda=0} =1.
\]

The representation $\rho^\vee_w$ of $\hat M^*=\GL_d(\C)\times \GL_d(\C)$ is given by $\Std\boxtimes \Std^\vee$,
where $\Std$ is the standard representation of $\GL_d$. Hence
\[
\rho^\vee_w\circ \phi_M = (\triv\boxtimes S_d) \otimes(\triv\boxtimes S_d)= \bigoplus_{i=0}^{d-1}\triv\boxtimes S_{2i+1}.
\]
From this, we compute
\begin{eqnarray*}
 r_P(w,\phi_{M,\lambda}) &=& \prod_{i=0}^{d-1} \frac{L(2\lambda,\triv\boxtimes S_{2i+1})}{\epsilon(2\lambda,\triv\boxtimes S_{2i+1})L(1+2\lambda,\triv\boxtimes S_{2i+1})}   \\
 &=& \prod_{i=0}^{d-1}  \frac{(1-q^{-2\lambda-i})^{-1}}{(-q^{\frac{1}{2}-2\lambda})^{2i}(1-q^{-1-2\lambda-i})^{-1}} ~=~ q^{(4\lambda-1)d(d-1)/2}\cdot \frac{1-q^{-2\lambda-d}}{1-q^{-2\lambda}}.
\end{eqnarray*}
Here, we use the formulas of local factors from \cite[Appendix A.1]{AGIKMS}. (The convention for the local $L$-factors and class field theory in \cite{AGIKMS} is opposite to that of \cite{KMSW}. The two instances of opposite signs cancel out in light of \cite[App.~F.1]{AGIKMS}.)
On the other hand, $\vol(\cO_D^\times)=\tfrac{q^d-1}{q^d}\vol(\cO_D)=\tfrac{q^d-1}{q^d} q^{-d(d-1)/2}$.
Therefore
  \[
r_P(w,\phi_{M,\lambda})^{-1}\cdot\frac{\vol(\cO_D^\times)}{1-q^{-2\lambda}} \bigg\rvert_{\lambda=0}
= \frac{1}{1-q^{-d}}\cdot \frac{q^d-1}{q^d}=1.
  \]  
This completes the proof.
\end{proof}

\section{Proof of LIR: the non-tempered case}\label{s:proof-non-tempered}

We treat the non-tempered case of Lemma~\ref{lem:reduction-to-maximal-Levi}.

\subsection{Reduction via globalization}

\begin{lem}\label{lem:reduction-nontempered}
In the non-tempered case, it suffices to prove LIR for the following data.
\begin{itemize}
\item $F$ is non-Archimedean,
\item $M=\GL_{1}(D)\times \GL_{1}(D)$,
\item $\psi_M=(\psi_M^\circ,\psi_M^\circ)$, where
$\psi^\circ_M\in \Psi_2(\GL_{d})$ is $D$-compatible and has the form $\rho\boxtimes S_1\boxtimes S_b$ (thus necessarily $d=b\dim\rho$),
\item $w\in W_{\psi_M}$ is the unique nontrivial element.
\end{itemize}
\end{lem}

\begin{proof}
The split case follows from Lemma~\ref{lem:easy-case}. For $d>1$,
the proof of Lemma~\ref{lem:reduction-tempered} applies with
case (ii) of Lemma~\ref{lem:globalization-parameter} in place of case (i).
\end{proof}

\subsection{Another explicit computation}

The proof of LIR is now reduced to the following proposition. A key idea in its proof is to exploit the compatibility between endoscopy and the Aubert--Zelevinsky duality, under which a representation of $D^\times$ corresponds to a tempered $L$-parameter and a non-tempered $A$-parameter (whose restriction to $\SL_2$ is nontrivial) at the same time.

\begin{prop}\label{prop:LIR-nontempered-special-case-1}
LIR is true in the case specified by Lemma \ref{lem:reduction-nontempered}.
\end{prop}

\begin{proof}
We put ourselves in the special case of Lemma \ref{lem:reduction-nontempered}. Let $\pi^\circ:=\pi_{\psi^\circ_M}\in \Irr(D^\times)$ correspond to the $A$-parameter $\psi^\circ_M=\rho \boxtimes S_1\boxtimes S_b$ via Proposition \ref{prop:local-JL}. We claim that the $L$-parameter of $\pi^\circ$ is $\phi^\circ_M:=\rho\boxtimes S_b$, the dual parameter of $\psi^\circ_M$. 

We prove the claim. Let $\pi^*:=\pi^*_{\psi^\circ_M}\in \Irr(\GL_d)$ be the image of $\psi^\circ_M$ under \eqref{eq:A-packet-GLn}.
Denote by $\hat \pi^\circ\in \Irr(D^\times)$ and $\hat \pi^*\in \Irr(\GL_d)$ their Aubert--Zelevinsky duals given by \cite{Zel80} for $\GL_d$ and \cite{Aubert} for general reductive groups. (Our hat corresponds to $(\cdot)^\#$ in \cite{Aubert}.) Since $D^\times$ is anisotropic modulo center, $\hat \pi^\circ=\pi^\circ$. On the other hand, $\hat \pi^{*}$ corresponds to the $A$-parameter $\phi^\circ_M=\rho\boxtimes S_b$ via \eqref{eq:A-packet-GLn} by \eqref{eq:Apacket-vs-Aubert}.
Now, Proposition \ref{prop:local-JL} for $\psi_M^\circ$ tells us that
\[
\Theta_{\pi^\circ}(g)= (-1)^{b-1} e(D^\times) \Theta_{\pi^{*}}(g^*), \qquad g\in D^\times ~\leftrightarrow~ g^*\in \GL_d(F).
\]
Applying Hiraga's theorem \cite[Thm.~1.5]{Hiraga} (with $G=D^\times$ and $H=\GL_d$ in his notation), since $D_{D^\times}(\pi^\circ)= \pi^\circ$, $D_{\GL_d}(\pi^{*})=(-1)^{d-b}\hat \pi^{*}$, and $(-1)^{a(D^\times)-a(\GL_d)}=(-1)^{1-d}$, we deduce that
\[
\Theta_{\hat \pi^\circ}(g)= e(D^\times) \Theta_{\hat \pi^{*}}(g^*), \qquad g\in D^\times ~\leftrightarrow~ g^*\in \GL_d(F).
\]
Comparing with Proposition \ref{prop:local-JL} for the parameter $\phi^\circ_M$, we conclude that $\hat \pi^\circ$ is the image of $\phi^\circ_M$ under the map of that proposition. That is, the claim is true.

We have already verified LIR for $G$, $M$, and $\phi_M=(\phi^\circ_M,\phi^\circ_M)$ in \S\ref{s:proof-tempered}. Therefore LIR for $G$, $M$, and $\psi_M$ is reduced to showing that the ratio of the normalized intertwining operators for $\phi_M$ and $\psi_M$ equals 1. Since they share the same induced representation and the same intertwining integral, it suffices to prove that
\begin{equation}\label{eq:ratios-equal}
    r_P(w,\psi_{M},\psi_F)/r_P(w,\phi_M,\psi_F)=\pi_M(\breve{w},\psi_M)/\pi_M(\breve{w},\phi_M),
\end{equation}
where 
the left hand side is understood as the value at $\lambda = 0$ of a meromorphic quotient, and 
the right hand side means the scalar ratio of the two operators. 

The right hand side equals $a_{\psi_M^\circ}/a_{\phi_M^\circ}\in\{\pm1\}$, as we see from Example \ref{ex:weyl-action} applied to each of $\psi_M^\circ$ and $\phi_M^\circ$. It follows from Lemma~\ref{lem:a-psi} that $a_{\psi_M^\circ}/a_{\phi_M^\circ}=(-1)^{b-1}$.
Therefore the proof is completed by the following lemma.
\end{proof}

\begin{lem}\label{lem:rP-ratio}
 In \eqref{eq:ratios-equal}, the left hand side equals $(-1)^{b-1}$.
\end{lem}

\begin{proof}
Recall that 
\[
r_P(w,\psi_{M,\lambda},\psi_F)=\frac{L(2\lambda,\rho^\vee_w\circ \phi_{\psi_{M}})}{\epsilon(2\lambda,\rho^\vee_w\circ \phi_{\psi_{M}})L(1+2\lambda,\rho^\vee_w\circ \phi_{\psi_{M}})}
\]
and
\[
\rho^\vee_w\circ \psi_{M}
= (\rho \boxtimes S_1 \boxtimes S_b)
\otimes (\rho^\vee \boxtimes S_1 \boxtimes S_b)
= \bigoplus_{i=0}^{b-1} \tau \boxtimes S_1 \boxtimes S_{2i+1},
\]
where we set $\tau = \rho \otimes \rho^\vee$.
Let $X(\tau)$ be the multiset over $\C/2\pi\sqrt{-1}\log(q)^{-1}\Z$ such that 
\[
\tau^{I_F} = \bigoplus_{\mu \in X(\tau)} |\cdot|^\mu, 
\]
where $I_F$ is the inertia subgroup.
Then 
\begin{align*}
&\frac{r_P(w,\psi_{M,\lambda},\psi_F)}{r_P(w,\phi_{M,\lambda},\psi_F)}   
\\&= \prod_{i=0}^{b-1} \frac{\epsilon(2\lambda, \tau \boxtimes S_{2i+1},\psi_F)}{\epsilon(2\lambda, \tau \boxtimes S_1 \boxtimes S_{2i+1},\psi_F)}
\frac{L(2\lambda, \tau \boxtimes S_1 \boxtimes S_{2i+1})}{L(2\lambda, \tau \boxtimes S_{2i+1})}
\frac{L(1+2\lambda, \tau \boxtimes S_{2i+1})}{L(1+2\lambda, \tau \boxtimes S_1 \boxtimes S_{2i+1})}
\\&=
\prod_{i=0}^{b-1} 
\det(-q^{\frac{1}{2}-2\lambda} \tau(\Frob)|\tau^{I_F})^{2i}
\frac{L(2\lambda-i,\tau)L(2\lambda-i+1,\tau) \cdots L(2\lambda+i-1,\tau)}
{L(1+2\lambda-i,\tau)L(1+2\lambda-i+1,\tau)\cdots L(1+2\lambda+i-1,\tau)}
\\&=
\prod_{i=0}^{b-1}  
\det(-q^{\frac{1}{2}-2\lambda} \tau(\Frob)|\tau^{I_F})^{2i}
\frac{L(2\lambda-i,\tau)}{L(2\lambda+i,\tau)}
\\&=
\prod_{i=0}^{b-1} 
\prod_{\mu \in X(\tau)}
(-q^{\frac{1}{2}-2\lambda-\mu})^{2i}
\frac{1-q^{-2\lambda-i-\mu}}{1-q^{-2\lambda+i-\mu}}.
\end{align*}
The factor for $i=0$ is $1$. Specializing at $\lambda=0$ gives
\begin{align*}
\frac{r_P(w,\psi_{M},\psi_F)}{r_P(w,\phi_{M},\psi_F)}   
&=
\prod_{i=1}^{b-1} 
\prod_{\mu \in X(\tau)}
q^{(1-2\mu)i}
\frac{1-q^{-i-\mu}}{1-q^{i-\mu}}.
\end{align*}
If we set 
\[
Q_i(\mu) = q^{(1-2\mu)i}
\frac{1-q^{-i-\mu}}{1-q^{i-\mu}}, 
\]
then 
\[
Q_i(\mu)Q_i(-\mu) = Q_i(\mu_0) = 1,
\quad Q_i(0) = -1, 
\]
where $\mu_0 \in \C/2\pi\sqrt{-1}\log (q)^{-1}\Z$ is such that $q^{-\mu_0} = -1$. 
Since $\tau$ is self-dual, $X(\tau)$ is invariant under $\mu \mapsto -\mu$. 
Hence 
\[
\prod_{\mu \in X(\tau)}Q_i(\mu) 
= (-1)^{\dim_\C \Hom_{W_F}(\C, \tau)} 
= (-1)^{\dim_\C \End_{W_F}(\rho)} 
= -1. 
\]
Therefore, we obtain that 
\[
\frac{r_P(w,\psi_{M},\psi_F)}{r_P(w,\phi_{M},\psi_F)} 
= \prod_{i=1}^{b-1}(-1) = (-1)^{b-1}, 
\]
as desired.
\end{proof}

\appendix
\section{Regularity and multiplicativity of Arthur-normalized intertwining operators}
\label{akms-reg:appendix}
In this appendix we prove Theorem~\ref{thm:regularity},  Lemma~\ref{lem:mult}, and the $\GL_m(D)$-case of \cite[Lem.~2.2.4]{KMSW} that was left unproven in that reference. 

Let $F$ be a local field of characteristic zero, $D$ a central
division algebra of degree $d$ over $F$, and
$G=\GL_m(D)$, $G^*=\GL_{dm}(F)$.
Write $\nu$ for $|\operatorname{Nrd}|_F$ on general linear groups
over $D$ and for $|\cdot|_F$ on $W_F$.
Fix a nontrivial unitary additive character $\psi_F$ of $F$.

Fix a Levi subgroup $M$ with
\[
 M=\prod_{i=1}^r\GL_{m_i}(D),\qquad \sum_i m_i=m,
\]
and a $D$-compatible Arthur parameter
$\psi_M=\boxtimes_i\psi_i$ with associated representation
$\pi_M=\boxtimes_i\pi_i$ as in Proposition~\ref{prop:local-JL}.
Put
\[
 \pi_{M,\lambda}=\boxtimes_i\pi_i\nu^{\lambda_i}
 \quad(\lambda\in\C^r),\qquad
 A_i=\phi_{\psi_i}=\phi_{\pi_{\psi_i}^*},\qquad
 B_i=\phi_{\pi_i}.
\]
Here $\phi_{\pi_i}$ is the $L$-parameter of $\pi_i$ under the local
Langlands correspondence for $\GL_{m_i}(D)$.

Write $\mathcal P(M)$ for the set of $F$-parabolic subgroups of $G$ with Levi component $M$. For $P_0=MN_0,P_1=MN_1\in\mathcal P(M)$, put
\[
 [J_{P_1\mid P_0}(\pi_{M,\lambda})f](g)
 =\int_{(N_0\cap N_1)\backslash N_1} f(ug)\,du
\]
in its region of convergence and extend meromorphically in $\lambda$.
We use the measures of \S\ref{ss:normalized-intertwining-operators}.
In \S\ref{akms-reg:sec:analytic} we normalize $J_{P_1\mid P_0}$
using the two collections $(A_i)$ and $(B_i)$, respectively.

\subsection{Comparison of local factors}
\label{akms-reg:sec:factors}
For a Langlands parameter $V$, set
\[
 \gamma(s,V,\psi_F)=\epsilon(s,V,\psi_F)
                    \frac{L(1-s,V^\vee)}{L(s,V)},
\]
\begin{equation}\label{akms-reg:eq:c}
 c(s,V,\psi_F)=\frac{L(s,V)}
 {\epsilon(s,V,\psi_F)L(1+s,V)},
\end{equation}
and
\begin{equation}\label{akms-reg:eq:k}
 r(s,V,\psi_F)=\epsilon(\tfrac12+s,V,\psi_F)c(s,V,\psi_F).
\end{equation}
Thus $r$ is the factor defining $r_{P'\mid P}$ in
\cite[\S2.2]{KMSW}. By \cite[(3.4.4), (3.4.7), (4.1.6)]{Tat79},
\begin{equation}\label{akms-reg:eq:opposite-k}
 r(s,V,\psi_F)\cdot r(-s,V^\vee,\psi_F)
 =\frac{\det V(-1)}
 {\gamma(s,V,\psi_F)\gamma(-s,V^\vee,\psi_F)}.
\end{equation}
We will also use
\begin{equation}\label{akms-reg:eq:k-positive}
 r(s,V,\psi_F)=b(V,\psi_F)\frac{L(s,V)}{L(1+s,V)},\qquad
 b(V,\psi_F)=b(V^\vee,\psi_F)>0,
\end{equation}
where $b(V,\psi_F)$ is independent of $s$.
For non-Archimedean $F$ and $\psi_F$ of conductor $\cO_F$,
$b(V,\psi_F)=q^{-a(V)/2}$, where $a(V)$ is the Artin conductor.
For the standard Archimedean characters, $b(V,\psi_F)=1$.
Replacing $\psi_F$ by $\psi_F(a\,\cdot)$ multiplies $b$ by
$|a|_F^{\dim V/2}$; see \cite[(3.4.4)--(3.4.7), (4.1.6)]{Tat79}.

\begin{lem}\label{akms-reg:lem:exact-comparison}
For all $i,j$,
\begin{equation}\label{akms-reg:eq:opposite-comparison}
\begin{aligned}
 &r(s,A_i\otimes A_j^\vee,\psi_F)\,
  r(-s,A_j\otimes A_i^\vee,\psi_F)\\
 &\qquad =
  r(s,B_i\otimes B_j^\vee,\psi_F)\,
  r(-s,B_j\otimes B_i^\vee,\psi_F).
\end{aligned}
\end{equation}
\end{lem}
\begin{proof}
Suppose first that $F$ is non-Archimedean. The Weil representation
in the Weil--Deligne parameter of an irreducible representation is
determined by its cuspidal support \cite[Prop.~4.3]{regCohen}.
The expansion of an irreducible class in standard module classes
preserves cuspidal support. Extended Jacquet--Langlands sends a
standard module to zero or transfers its essentially square-integrable
factors \cite[\S4.1]{regCohen}, preserving their Langlands parameters.
Thus $\operatorname{LJ}([\pi_{\psi_i}^*])=\pm[\pi_i]$ implies that
$A_i$ and $B_i$ have isomorphic underlying Weil representations.
The determinants and gamma factors of the tensor products therefore
agree \cite[Cor.~4.5(1)]{HM17}. Now use \eqref{akms-reg:eq:opposite-k}.

For split Archimedean groups, $A_i=B_i$. For $(F,D)=(\R,\mathbb H)$,
write
\[
 U=\nu^{1/2}\oplus\nu^{-1/2},\qquad
 \rho_m=\operatorname{Ind}_{W_\C}^{W_\R}(z/|z|)^m\ (m\geq1),
 \qquad \rho_0=1\oplus\sgn .
\]
By \cite[\S11.3, Cor.~13.7, Thm.~13.8]{BadRen}, the constituents
$\rho_m\nu^{it}\boxtimes S_b$ of $\psi_i$ contribute equally to
$A_i$ and $B_i$. A constituent $\chi\boxtimes S_{2a}$ with
$\dim\chi=1$ gives $\chi\circ\det$ and $\chi\circ\operatorname{Nrd}$,
respectively. In their Langlands parameters this replaces
$\chi\nu^{a-1-2j}U$ by $\chi\nu^{a-1-2j}\rho_1$ for
$0\leq j<a$, using $\rho_1\otimes\sgn\simeq\rho_1$.

For an irreducible $W_\R$-representation $T$, the formulas for
$L$-factors \cite[(3.2.4)--(3.2.5)]{Tat79} give
\[
 \frac{L(s,U\otimes T)}{L(s,\rho_1\otimes T)}
 =
 \begin{cases}
  \dfrac{2\pi}{s+t-\tfrac12}, & T=\nu^t,\\[3pt]
  1, & T=\sgn\,\nu^t\text{ or }\rho_m\nu^t\ (m\geq1).
 \end{cases}
\]
Use $\Gamma_\R(z+2)=\frac{z}{2\pi}\Gamma_\R(z)$ and
$\rho_1\otimes\rho_m=\rho_{m+1}\oplus\rho_{m-1}$.
The two tensors have the same dimension, so their $b$-factors
in \eqref{akms-reg:eq:k-positive} agree. The ratio of $r$-factors
is $(s+t+\tfrac12)/(s+t-\tfrac12)$ for $T=\nu^t$ and $1$ otherwise.
Replacing $(s,T)$ by $(-s,T^\vee)$ inverts this ratio. Since
$U$ and $\rho_1$ are self-dual, each replacement in
either tensor factor leaves the product unchanged. Twists are
absorbed into $T$, and multiplicativity proves the identity.
\end{proof}

\subsection{Regularity and transitivity of the normalized intertwining operators}
\label{akms-reg:sec:analytic}
For $P_0,P_1\in\mathcal P(M)$, define $i\prec_{P_k}j$ when
$N_{P_k}$ contains the matrix root subgroup in position $(i,j)$.
Put
\[
 X(P_1,P_0)=\{(i,j):i\prec_{P_0}j,\ j\prec_{P_1}i\}.
\]

For $C=A,B$, define
\begin{align*}
 r^C_{P_1\mid P_0}(\lambda,\psi_F)
 &=\prod_{(i,j)\in X(P_1,P_0)}
 r(\lambda_i-\lambda_j,C_i\otimes C_j^\vee,\psi_F),\\
 R^C_{P_1\mid P_0}(\lambda,\psi_F)
 &=(r^C_{P_1\mid P_0}(\lambda,\psi_F))^{-1}
 J_{P_1\mid P_0}(\pi_{M,\lambda}).
\end{align*}

\begin{lem}\label{akms-reg:lem:exact-rankone}
Let $P$ be the upper-block-triangular parabolic in $\GL_m(D)$, where $m=a+b$, 
with Levi $M=\GL_a(D)\times\GL_b(D)$, and let $\overline P$ be its
opposite. For irreducible $\pi_1,\pi_2$ with Langlands parameters
$B_1,B_2$, put $V=B_1\otimes B_2^\vee$. Then
\begin{equation}\label{akms-reg:eq:exact-rankone}
 J_{P\mid\overline P}(\pi_1\nu^s\boxtimes\pi_2)
 J_{\overline P\mid P}(\pi_1\nu^s\boxtimes\pi_2)
 =r(s,V,\psi_F)\cdot r(-s,V^\vee,\psi_F)\cdot \operatorname{id}_{\mathcal I_P(\pi_1\nu^s\boxtimes\pi_2)}.
\end{equation}
\end{lem}
\begin{proof}
First let $\pi_1,\pi_2$ be unitary square-integrable.
For $F=\R$, use \cite[Lem.~2.2.1 and its proof]{KMSW}; for
$F=\C$, use \cite[Thm.~2.1, \S3 (3.2)]{ArtIOR1} for
$\operatorname{Res}_{\C/\R}\GL_n$. The identity for all $s$ follows
by meromorphic continuation from unitary twists.

Suppose $F$ is non-Archimedean and $\psi_F$ has conductor $\cO_F$.
Compare $\GL_{a+b}(D)$ with $\GL_{d(a+b)}(F)$, matching the
inducing representations by Jacquet--Langlands. With
$\operatorname{vol}(\cO_D)=\operatorname{vol}(\cO_F)=1$, write
$J_{P\mid\overline P}J_{\overline P\mid P}=j_{\GL_m(D)}^0\operatorname{id}$,
and similarly $j_{\GL_{dm}(F)}^0$ for the split group.
Put
\[
 C_r(Q)=\prod_{v=1}^r(1-Q^{-v}),\qquad
 g_{\GL_{m}(D)}=\frac{C_{a+b}(q^d)}{C_a(q^d)C_b(q^d)},\qquad
 g_{\GL_{dm}(F)}=\frac{C_{d(a+b)}(q)}{C_{da}(q)C_{db}(q)}.
\]
By \cite[\S2 (1), Thms.~3.1, 7.2]{AubertPlymen}, with the same
congruence-subgroup calculation over $\cO_D$,
\begin{gather*}
 \Lambda_d(r)=q^{-d(d-1)r^2/2}\frac{C_r(q^d)}{C_{dr}(q)},\\
 \frac{j_{\GL_m(D)}^0}{j_{\GL_{dm}(F)}^0}
 =\frac{g_{\GL_{m}(D)}}{g_{\GL_{dm}(F)}}\cdot
   \frac{\Lambda_d(a)\Lambda_d(b)}{\Lambda_d(a+b)}
 =q^{d(d-1)ab}.
\end{gather*}
Here formal degrees transfer by $\Lambda_d(r)$.
The canonical orbit measures of \cite[\S2]{AubertPlymen} agree,
since $\operatorname{Nrd}(tI_r)=t^{dr}=\det(tI_{dr})$ and the
unramified stabilizers agree by \cite[\S7 (48)]{AubertPlymen}.
For the self-dual measures,
$\operatorname{vol}(\cO_D)=q^{-d(d-1)/2}$
(Lemma~\ref{lem:lift-and-vol}(ii)); each of the two integrations has
$ab$ coordinates over $D$. Hence
\[
 j_{\GL_m(D)}=q^{-d(d-1)ab}\cdot j_{\GL_m(D)}^0=j_{\GL_{dm}(F)}^0=j_{\GL_{dm}(F)}.
\]
On the split side, \cite[\S3 (4), \S4 (7)--(9)]{AubertPlymen}
gives \eqref{akms-reg:eq:exact-rankone}, since $c(G\mid M)=1$ and (assuming the conductor of $\psi_F$ is $\mathcal O_F$)
\[
 r(s,V,\psi_F)=q^{-a(V)/2}\frac{L(s,V)}{L(1+s,V)},
\]
where $a(V)$ is the Artin conductor. See also
\cite[(25.1), Lem.~25.2, p.~169]{Art05}.
Changing $\psi_F$ scales both sides equally, so the condition on the conductor of $\psi_F$ can be dropped.

The case of unitary square-integrable $\pi_i$ is thus complete. Since both sides of the identity behave well under twisting by characters, the case of not necessarily unitary but essentially square-integrable $\pi_i$ follows immediately. For arbitrary $\pi_i$, choose standard modules
\[
 S_i=\delta_{i,1}\nu^{x_{i,1}}\times\cdots\times
     \delta_{i,t_i}\nu^{x_{i,t_i}}\twoheadrightarrow\pi_i,
 \qquad
 B_i=\bigoplus_u\phi_{\delta_{i,u}}\nu^{x_{i,u}}
\]
with unitary square-integrable $\delta_{i,u}$, as in
\cite[Thm.~3.3]{BadRen}.
The integrals defining $J$ commute with these quotient maps,
first where they converge and then meromorphically.
By \cite[\S1 (J2)--(J3)]{ArtIOR1}, their product on
$S_1\nu^s\times S_2$ is the product of
\eqref{akms-reg:eq:exact-rankone} for the pairs
$(\delta_{1,u},\delta_{2,v})$ at $s+x_{1,u}-x_{2,v}$.
Multiplicativity of Artin factors gives
$r(s,V)r(-s,V^\vee)$. The identity descends to
$\pi_1\nu^s\times\pi_2$.
\end{proof}

\begin{prop}\label{akms-reg:prop:composition}
For $C=A,B$, the operators
$R^C_{P_1\mid P_0}(\lambda)$ are regular and invertible
at $\lambda=0$. For $P_0,P_1,P_2\in\mathcal P(M)$, they satisfy
\begin{equation}\label{akms-reg:eq:composition}
 R^C_{P_2\mid P_1}(\lambda)\circ 
 R^C_{P_1\mid P_0}(\lambda)
 =R^C_{P_2\mid P_0}(\lambda).
\end{equation}
\end{prop}
\begin{proof}
In the two-block case put $\lambda=(s,0)$ up to a common twist and
$V=C_1\otimes C_2^\vee$. Lemma~\ref{akms-reg:lem:exact-rankone} and
\eqref{akms-reg:eq:opposite-comparison} give
\begin{equation}\label{akms-reg:eq:inverse}
 R^C_{P_0\mid\overline P_0}(s,0)
 R^C_{\overline P_0\mid P_0}(s,0)=1.
\end{equation}
Since the $C_i$ are parameters of unitary representations, the real
exponents in their tempered decompositions occur in opposite pairs
with the same tempered parameters. Thus $V^\vee\simeq\overline V$,
and \eqref{akms-reg:eq:k-positive} gives
$r(-s,V^\vee)=\overline{r(s,V)}$ for $s\in i\R$. By
\cite[\S1 (J4)]{ArtIOR1}, the operators in
\eqref{akms-reg:eq:inverse} are adjoints on $i\R$ wherever defined.
They are therefore unitary there. Their meromorphic matrix entries
on each finite-dimensional isotypic component for a fixed maximal
compact subgroup are bounded along $i\R$, hence regular at zero. The inverse
identity gives invertibility at zero.

For general $P_0,P_1$ use a reduced sequence of adjacent swaps of
blocks. By \cite[\S1 (J2)--(J3)]{ArtIOR1}, the integrals factor into
the corresponding two-block operators. Each pair in $X(P_1,P_0)$
occurs once, so the normalizing factors multiply to
$r^C_{P_1\mid P_0}$. This gives regularity and invertibility at zero.
The reduced products depend only on their endpoints, and
\eqref{akms-reg:eq:inverse} allows cancellation of a swap followed
by its reverse. The relations in the symmetric group give
\eqref{akms-reg:eq:composition}.
\end{proof}

\begin{proof}[Proof of Theorem~\ref{thm:regularity}]
Put $Q=\breve w^{-1}P\breve w$ and
$[\ell(\breve w)f](g)=f(\breve w^{-1}g)$.
The change of variables $u\mapsto\breve w^{-1}u\breve w$ identifies
the integration measures, so
$J_P(\breve w,\pi_{M,\lambda})=\ell(\breve w)J_{Q\mid P}(\pi_{M,\lambda})$.
By \eqref{akms-reg:eq:k},
\[
 R_P(\breve w,\psi_{M,\lambda},\psi_F)
 =
 \left(\prod_{(i,j)\in X(Q,P)}
   \epsilon(\tfrac12+\lambda_i-\lambda_j,
             A_i\otimes A_j^\vee,\psi_F)\right)
 \ell(\breve w)R^A_{Q\mid P}(\lambda,\psi_F).
\]
The product is holomorphic and nowhere zero. The assertion follows
from Proposition~\ref{akms-reg:prop:composition}.
\end{proof}

\subsection{Multiplicativity of the normalized self-intertwining operator}
\label{akms-reg:sec:multiplicativity}

Write $\sigma(w)=[n(w)]$ for the section of \eqref{eq:splnpsi},
where $n(w)=\widehat{\tilde w}$. Let $\psi_t$ range over the
distinct block parameters and let $w_t$ be the restriction of $w$
to the blocks with parameter $\psi_t$. Use identical representation
models on these blocks and put
$T_w(\bigotimes_i v_i)=\bigotimes_i v_{w^{-1}(i)}$.

\begin{lem}\label{akms-reg:lem:pi-composition}
For $w\in W_{\psi_M}$,
\begin{equation}\label{akms-reg:eq:signed-permutation}
 \pi_M(\breve w,\psi_M)=\kappa(w)T_w,\qquad
 \kappa(w)=\prod_t \sgn (w_t)^{(1-a_{\psi_t})/2}.
\end{equation}
\end{lem}
\begin{proof}
Since $\sigma$ is a homomorphism, $n(w^{-1})^{-1}$ represents
$\sigma(w)$. For this representative the element
$s=u^{-1}n(w^{-1})^{-1}$ in \cite[\S2.4.3]{KMSW} is $1$.

On a cycle of length $\ell$ with common parameter $\psi_0$ on a
factor $H$, the twisted character identity reduces to the ordinary
character identity for $H$ with endoscopic element $s_{\psi_0}^{\ell}$
\cite[\S1.1.3, Cor.~1.1.6, Lem.~1.1.7]{KMSW}. For
$\sigma_0=\pi_{\psi_0}$ its coefficient is
\[
 e(H)\langle s_{\psi_0}^{\ell+1},\sigma_0\rangle
 =e(H)a_{\psi_0}^{\ell-1}
\]
by \cite[Thm.~1.6.1(4), (1.6.1)]{KMSW} and $s_{\psi_0}^2=1$.
The product of the $e(H)$ over the cycles is $e(M^\theta)$.
Thus \eqref{eq:check-w-characterization} gives the factor
$a_{\psi_0}^{\ell-1}$ on each cycle. Multiplying over the cycles
gives \eqref{akms-reg:eq:signed-permutation}.
\end{proof}

\begin{cor}\label{akms-reg:cor:compound-multiplicativity}
For $w_1,w_2\in W_{\psi_M}$,
\[
 R_P(w_1w_2,\psi_M)=R_P(w_1,\psi_M)R_P(w_2,\psi_M).
\]
\end{cor}
\begin{proof}
Both factors in \eqref{akms-reg:eq:signed-permutation} are
multiplicative in $w$. By \eqref{eq:RP(w)} and
\eqref{akms-reg:eq:k}, the first arrow is $R^A_{w^{-1}Pw\mid P}$
followed by the normalized translation of
\cite[\S2.3]{KMSW}. Apply the calculation in
\cite[Lem.~2.5.3]{KMSW} using \eqref{akms-reg:eq:composition},
\cite[Lem.~2.2.5]{KMSW}, and
\cite[Lem.~2.3.1, Fact~2.3.2]{KMSW}.
\end{proof}

\section{Compatibility with parabolic induction}\label{akms-ind:appendix}

Write $\psi_M=(\psi_i)$. Choose ordered decompositions
$\psi_i=\bigoplus_\alpha\eta_{i,\alpha}$ into irreducibles, using the
same list for isomorphic $\psi_i$. Each constituent is $D$-compatible.
Let $Q\subset M$ be the corresponding standard parabolic with Levi
$M_0$, and put $P_0=QN_P$ and $\psi_0=(\eta_{i,\alpha})_{i,\alpha}$.
For $w\in W_{\psi_M}$ define $w_0\in W_{\psi_0}$ by
$w_0(i,\alpha)=(w(i),\alpha)$.

\begin{lem}\label{akms-ind:lem:refinement}
Under normalized induction in stages,
\[
 R_P(w,\psi_M)=R_{P_0}(w_0,\psi_0).
\]
\end{lem}
\begin{proof}
Put $\tau=\pi_{\psi_0}$. Proposition~\ref{prop:local-JL} and the
character formula for normalized induction give
$\pi_M=\cI_Q^M(\tau)$ and $a_{\psi_i}=\prod_\alpha a_{\eta_{i,\alpha}}$,
since $e(\GL_h(D))=(-1)^{(d-1)h}$. Under induction in stages, the
isomorphism is given by $(\Omega F)(g)=F(g)(1)$.

By Lemma~\ref{lem:independent}, we may use the isomorphism $\xi$ of
\S\ref{ss:notation} for both pairs, with the same cocycle for $D^\times$
on each diagonal $d\times d$ block. Both $w$ and $w_0$ send coordinate
$k$ of block $i$ of $M^*$ to coordinate $k$ of block $w(i)$. Thus
their lifts agree, and $n:=\breve w=\breve w_0$ normalizes $Q$.
Since the order of the coordinates within each block of $M^*$ is
preserved, the inversion roots for $P$ and $P_0$ coincide. Hence
\[
 N_{P_0}\cap nN_{\overline P_0}n^{-1}
 =N_P\cap nN_{\overline P}n^{-1},\qquad
 \rho_{w_0}\simeq\rho_w|_{{}^LM_0^*}.
\]
The integrals thus have the same domain and measure, and their
normalizing factors agree. Twist each $\eta_{i,\alpha}$ and the
corresponding representation by $|\cdot|^{\lambda_{i,\alpha}}$.
The map $\Omega$ identifies the normalized integrals in a common
chamber of convergence, so they agree by meromorphic continuation.
The twists of $M$ correspond to $\lambda_{i,\alpha}=\lambda_i$.
Specializing at zero by Theorem~\ref{thm:regularity} identifies the
first arrows in \eqref{eq:RP(w)} for $P$ and $P_0$.

On a pure tensor of induced functions,
$(T_w\varphi)(1)=T_{w_0}(\varphi(1))$. Each cycle of $w$ of length
$\ell$ with parameter $\psi_i$ gives one cycle of $w_0$ of length
$\ell$ for each $\eta_{i,\alpha}$. The corresponding factors in
\eqref{akms-reg:eq:signed-permutation} agree since
$a_{\psi_i}^{\ell-1}=\prod_\alpha a_{\eta_{i,\alpha}}^{\ell-1}$.
Thus $\Omega$ also identifies the second arrows in \eqref{eq:RP(w)}.
\end{proof}

\bibliographystyle{alpha}
\bibliography{bib}

\end{document}